\documentclass[a4paper,10pt]{article}
\usepackage{amssymb}
\usepackage{latexsym}
\usepackage{amsmath}
\usepackage{amsthm}
\usepackage{amsfonts}
\usepackage{amssymb}
\usepackage{amsmath,amscd}
\usepackage{graphicx}

\newtheorem{Th}{Theorem}

\newtheorem{Lm}{Lemma}
\newtheorem{Lma}{Lemma}[section]

\newcommand{\be}{\begin{equation}}
\newcommand{\ee}{\end{equation}}
\newcommand{\bes}{\begin{equation*}}
\newcommand{\ees}{\end{equation*}}
\newcommand{\Di}{\mathbb{D}}
\newcommand{\R}{\mathbb{R}}
\newcommand{\N}{\mathbb{N}}

\newcommand{\reset}{\setcounter{equation}{0}\setcounter{Th}{0}\setcounter{Prop}{0}\setcounter{Co}{0}
\setcounter{Lm}{0}\setcounter{Rm}{0}}

\def\Xint#1{\mathchoice
{\XXint\displaystyle\textstyle{#1}}%
{\XXint\textstyle\scriptstyle{#1}}%
{\XXint\scriptstyle\scriptscriptstyle{#1}}%
{\XXint\scriptscriptstyle\scriptscriptstyle{#1}}%
\!\int}
\def\XXint#1#2#3{{\setbox0=\hbox{$#1{#2#3}{\int}$}
\vcenter{\hbox{$#2#3$}}\kern-.5\wd0}}
\def\dashint{\Xint-}

\def\ti{\tilde}
\def\lf{\left}
\def\rg{\right}

\def\al{\alpha}
\def\la{\lambda}
\def\e{\varepsilon}
\def\ep{\varepsilon}

\def\ds{\displaystyle}
\def\ov{\overline}
\def\Om{\Omega}
\def\om{\omega}
\def\p{\partial}

\makeatletter
\newcommand{\tpitchfork}{%
  \vbox{
    \baselineskip\z@skip
    \lineskip-.52ex
    \lineskiplimit\maxdimen
    \m@th
   \ialign{##\crcr\hidewidth\smash{$- $}\hidewidth\crcr$\pitchfork$\crcr}
  }%
}
\makeatother

\begin{document}

\title{Conservation Laws for $p $-Harmonic Systems with Antisymmetric Potentials and Applications}

\author{Francesca Da Lio 
 \ and Tristan Rivi\`ere\footnote{Department of Mathematics, ETH Zentrum,
CH-8093 Z\"urich, Switzerland.}}

\maketitle

\medskip

{\bf Abstract :} {\it  We prove that $p $-harmonic systems with antisymmetric potentials of the form
\[
-\,\mbox{div}\lf((1+|\nabla u|^2)^{\frac{p}{2}-1}\,\nabla u\rg)=(1+|\nabla u|^2)^{\frac{p}{2}-1}\,\Om\cdot \nabla u\ ,
\]
($\Om$ is antisymmetric) can be written in divergence form as a conservation law
\[
-\mbox{div}\lf((1+|\nabla u|^2)^{\frac{p}{2}-1}\,A\,\nabla u\rg)=\nabla^\perp B\cdot\nabla u\ .
\]
This extends to the $p $-harmonic framework the original work of the second author for $p=2$ (see \cite{Riv}).
We give applications of the existence of this divergence structure in the analysis $p\rightarrow 2$.}
\section{Introduction}

In \cite{Riv} the second author studied Linear Elliptic Systems in two dimensions with antisymmetric potentials of the form
\be
\label{0.1}
-\Delta u=\Om\cdot \nabla u\ ,
\ee
where $u\in W^{1,2}(\Di^2,{\R}^m)$ with $m\ge 2$ and $\Om\in L^2(\Di^2,so(m)\otimes {\R}^2)$ where $\Di^2$ denotes the two dimensional disc and $so(m)$ the space of $m\times m$
antisymmetric matrices. One of the motivations for considering these systems is that it is in particular satisfied by critical points to arbitrary conformally invariant coercive Lagrangians of maps with quadratic growth (harmonic maps into manifolds, prescribed mean curvature surfaces in euclidian and non euclidian spaces). While the system (\ref{0.1}) looks at a first glance critical with respect to the $W^{1,2} $-norm ($\Om\cdot\nabla u\in L^1(\Di^2,{\R}^m)$), the main contribution in \cite{Riv} was to stress the importance of the anti-symmetry of $\Om$ for being responsible of a-priori unexpected sub-critical behaviours under small $L^2$ norm of $\Om$ assumption (regularity, pre-compactness below arbitrary energy levels...etc)

The main idea for proving these sub-critical behaviour under small $L^2$ norm of $\Om$ assumption was the discovery of conservation laws for solutions to (\ref{0.1}) and the proof of the existence of $A\in W^{1,2}(\Di^2,Gl_m({\R}))$ together with $B\in W^{1,2}(\Di^2,M_m({\R}))$ (where $Gl_m({R})$ is denoting the group of invertible $m\times m$ real matrices and $M_m({\R})$ the space of square $m\times m$ real matrices) such that
\be
\label{0.2}
-\nabla A+A\,\Om=\nabla^\perp B\ ,
\ee
where $\nabla^\perp\cdot:=(-\p_{x_2}\cdot,\p_{x_1}\cdot)$ with the estimates
\be
\label{0.3}
\lf\|\mbox{dist}(A,SO(m))\rg\|^2_{L^\infty(\Di^2)}+\int_{\Di^2} |\nabla A|^2\ dx^2\le C\, \int_{\Di^2} |\Om|^2\ dx^2\ ,
\ee
where $C>0$ is a constant depending only on $m$. The existence of the pair $(A,B)$ permits to rewrite the system in conservative form
\be
\label{0.4}
-\mbox{div}(A\,\nabla u)=\nabla^\perp B\cdot\nabla u\ ,
\ee
from which the sub-critical behaviours follow thanks in particular to Wente type estimates (see \cite{Hel}). In \cite{Riv-1} the second author has been asking whether such an approach
could be available for considering the $p $-harmonic counterpart of (\ref{0.4}) in particular for solutions to systems of the form
\be
\label{0.5}
-\,\mbox{div}\lf((1+|\nabla u|^2)^{\frac{p}{2}-1}\,\nabla u\rg)=(1+|\nabla u|^2)^{\frac{p}{2}-1}\,\Om\cdot \nabla u\ ,
\ee
 for $p\ne 2$ and where $\Om$ is again antisymmetric. This is the purpose of the present work to present a contribution in this direction

\medskip

Let $p\ge 2$. We first consider critical points to the $p $- energy
\be
\label{IV.1}
E_p(u):=\int_{\Di^2}(1+|\nabla u|^2)^{\frac{p}{2}}\ dx^2
\ee
among maps taking values into the sphere $S^{m-1}\hookrightarrow {\R}^m$ that is in the space
\be
\label{IV.2}
W^{1,p}(\Di^2,S^{m-1}):=\lf\{  u\in W^{1,p}(\Di^2,{\R}^m)  \ ;\quad |u|^2=1\quad \mbox{ a.e. }\rg\}\ .
\ee
It is well known that such a map satisfies the Euler Lagrange equation
\be
\label{IV.3}
-\,\mbox{div}\lf((1+|\nabla u|^2)^{\frac{p}{2}-1}\,\nabla u\rg)= (1+|\nabla u|^2)^{\frac{p}{2}-1}\ u\ |\nabla u|^2 \quad\mbox{ in }{\mathcal D}'(\Di^2)\ .
\ee
This gives in particular
\be
\label{IV.4}
\mbox{div}\lf((1+|\nabla u|^2)^{\frac{p}{2}-1}\,u\wedge\nabla u\rg)= 0\ ,
\ee
which in coordinates gives
\be
\label{IV.5}
\forall\ i,j\in \{1\cdots m\}\quad\quad \mbox{div}\lf((1+|\nabla u|^2)^{\frac{p}{2}-1}\,\lf[u^i\,\nabla u^j-u^j\,\nabla u^i\rg]\rg)= 0\ .
\ee
We observe that for any $i\in\{1\cdots m\}$ there holds
\be
\label{IV.6}
\begin{array}{rl}
\ds(1+|\nabla u|^2)^{\frac{p}{2}-1}\ u^i\ |\nabla u|^2 &\ds=\sum_{j=1}^m(1+|\nabla u|^2)^{\frac{p}{2}-1}\ u^i\ \nabla u^j\cdot \nabla u^j \\[5mm]
 &\ds =\sum_{j=1}^m(1+|\nabla u|^2)^{\frac{p}{2}-1}\ \lf(u^i\ \nabla u^j-u^j\,\nabla u^i\rg)\cdot \nabla u^j 
\end{array}
\ee
where we have used
\be
\label{IV.7}
\sum_{j=1}^m u^j\cdot \nabla u^j =0\ .
\ee
Introducing
\be
\label{IV.8}
\nabla^\perp B^{ij}:=(1+|\nabla u|^2)^{\frac{p}{2}-1}\ \lf(u^i\ \nabla u^j-u^j\,\nabla u^i\rg)\ ,
\ee
we rewrite the Euler Lagrange Equation (\ref{IV.3}) as follows
\be
\label{IV.9}
-\,\mbox{div}\lf((1+|\nabla u|^2)^{\frac{p}{2}-1}\,\nabla u\rg)=\nabla^\perp B\cdot\nabla u\quad\mbox{ in }{\mathcal D}'(\Di^2)\ .
\ee
These observation have been used in particular in \cite{Strz} for proving the regularity of $n $-harmonic maps into the spheres.

\medskip

Modifying the target and considering critical points of $E_p$ but for maps taking values into an arbitrary sub-manifold of ${\R}^m$ leads to an Euler-Lagrange equation
of the form
\be
\label{0.6}
-\,\mbox{div}\lf((1+|\nabla u|^2)^{\frac{p}{2}-1}\,\nabla u\rg)=(1+|\nabla u|^2)^{\frac{p}{2}-1}\,\Om\cdot \nabla u\ .
\ee
The main result of the present work is to prove that the above conservation law (\ref{IV.9}) exist for solutions to general systems of the form (\ref{0.6}) (including in particular $p$-harmonic maps into general targets). Precisely we prove the following theorem.
\begin{Th}
\label{th-I.1}
Let $p\ge 2$. There exists $\sigma>0$ depending only on $m\in {\N}$ such that for any  $\Om\in L^2(\Di^2,{\R}^2\otimes so(m))$ and $u\in W^{1,p}(\Di^2,{\R}^m)$ satisfying
\be
\label{I.1}
-\,\mbox{div}\lf((1+|\nabla u|^2)^{\frac{p}{2}-1}\,\nabla u\rg)=(1+|\nabla u|^2)^{\frac{p}{2}-1}\,\Om\cdot \nabla u\ ,
\ee
if
 \be
\label{I.2}
\int_{\Di^2}(1+|\nabla u|^2)^{\frac{p}{2}-1}\ |\Om|^2\ dx^2<\sigma\ ,
\ee
then there exist $A\in L^\infty\cap W^{1,2}(\Di^2,Gl_m({\R}))$ and $B\in W^{1,p/(p-1)}(\Di^2,M_m({\R}))$ such that
\be
\label{I.3}
\lf\|\mbox{dist}(A,SO(m))\rg\|^2_{L^\infty(\Di^2)}+\int_{\Di^2}(1+|\nabla u|^2)^{\frac{p}{2}-1}\,|\nabla A|^2\ dx^2\le C\, \int_{\Di^2}(1+|\nabla u|^2)^{\frac{p}{2}-1}\ |\Om|^2\ dx^2
\ee
such that
\be
\label{I.4}
\int_{\Di^2}\frac{|\nabla B|^2}{(1+|\nabla u|^2)^{\frac{p}{2}-1}\,}\ dx^2\le C\, \int_{\Di^2}(1+|\nabla u|^2)^{\frac{p}{2}-1}\ |\Om|^2\ dx^2
\ee
and
\be
\label{I.5}
-\mbox{div}\lf((1+|\nabla u|^2)^{\frac{p}{2}-1}\,A\,\nabla u\rg)=\nabla^\perp B\cdot\nabla u\ .
\ee
where the constant $C$ above is not depending on $p$. \hfill $\Box$
\end{Th}
The Lagrangian $E_p$ for $p> 2$ in 2 dimension has first been introduced by J.Sacks and K.Uhlenbeck in \cite{SaU} in order to give a sub-critical relaxation of the homotopic Plateau 
Problem in arbitrary closed manifold. As $p$ converges to $2$ critical points to $E_p$ tend to form a bubble tree of conformal harmonic maps connected by ''necks'' . 

One of the main contribution of \cite{LR} was to show the importance of conservation laws of the form (\ref{0.4}) in order to give a precise description of the sequences
of maps solving (\ref{0.1}) in the neck regions connecting the bubbles. There are three main questions posed by this bubble formation
\begin{itemize}
\item[ ]{i)} Does energy dissipate in the neck regions ?
\item[ ]{ii)} Has the neck region some asymptotic non zero length ?
\item[ ]{iii)} Is there some part of the ``Morse index'' being dissipated in the neck region ?
\end{itemize}
The negative answer  to the first question is called {\em energy quantization} or  {\em no-neck energy} phenomenon. The negative answer  to the second question is called {\em necklessness} property
while the negative  answer to the third question is called {\em Morse index stability}. Since the pioneered work of Sacks and Uhlenbeck there has been numerous contributions aiming at giving answers to these  questions  (\cite{DK}, \cite{Lam}, \cite{LCM}, \cite{LMM}...)

In order to consider these three questions,  the analysis developed by P.Laurain and the second author in \cite{LR}, which was based on an idea that first appeared  in \cite{LiRi}, and  which consists in combining  dualities in interpolation Lorentz spaces with the existence of conservation laws, has shown to be very efficient for critical points to conformally invariant Lagrangians of mappings into manifolds\footnote{Since then, this powerful strategy of combining the $L^{p,\infty}-L^{p',1}$ duality with conservation laws  for studying the no neck property, the neckless property and the Morse index stability has been extensively used in different contexts beyond the framework of conformally invariant lagrangians of mappings into manifolds such as Yang-Mills fields \cite{Riv-0}, \cite{Gau} and Willmore surfaces \cite{BR0} and \cite{MiRi}.}.

In order to  illustrate possible use of the conservation laws (\ref{I.5}) to study the three questions i), ii) and iii) above in the context of Sacks-Uhlenbeck relaxation,   we prove the following result which is giving sufficient condition for controlling the length of the necks at the limit.
\begin{Th}
\label{th-IV.2-ann}
For any $p> 2$ the following holds. Let $0<\delta<1$ and denote
\be
\label{II.an-1}
\sqrt{p-2}\,\log\delta^{-1}:=M\ .
\ee
Let $\Om\in L^2(\Di^2,{\R}^2\otimes so(m))$ and $u\in W^{1,p}(\Di^2,{\R}^m)$ satisfying
\be
\label{II-r.1-ann}
-\,\mbox{div}\lf((1+|\nabla u|^2)^{\frac{p}{2}-1}\,\nabla u\rg)=(1+|\nabla u|^2)^{\frac{p}{2}-1}\,\Om\cdot \nabla u\ \mbox{ in }{\mathcal D}'(\Di^2)\ .
\ee
Denote
 \be
\label{II-r.2-ann}
\int_{B_1(0)\setminus B_\delta(0)}(1+|\nabla u|^2)^{\frac{p}{2}-1}\ |\Om|^2\ dx^2=Q\ .
\ee
and
\be
\label{II.r.3-ann}
\lf\|\max\{|x|,\delta\}\ |\nabla u|\rg\|_{L^\infty(B_1(0))}= L
\ee
then for ant $t<1$ we have
\be
\label{II.r.4-ann}
\lf\|u(x)-u(y)\rg\|_{L^\infty(B_t(0)\times B_t(0))}\le C(t,M,L,\|\nabla u\|_p,Q)\ .
\ee
\hfill $\Box$
\end{Th}
The conditions 
\be
\label{II.r7-ann}
\limsup_{p\rightarrow 2} \sqrt{p-2}\,\log\delta^{-1}<+\infty
\ee
and (\ref{II.r.3-ann}) for ensuring the neck length control were already known (see \cite{LiWa-1}) in the special case of $p $-harmonic maps
into smooth (at least $C^3$) manifolds. The present result is extending this fact to a much wider case of PDE's including $p $-harmonic relaxation of arbitrary conformally invariant Lagrangians of maps. Observe that the Sacks Uhlenbeck $\epsilon $-regularity for $p $-harmonic maps (independent of $p$ as $p$ tends to 2 from above) is implying (\ref{II.r.3-ann})  together with the condition
\be
\label{II.r5-ann}
\limsup_{p\rightarrow 2} \,(p-2)\,\log\delta^{-1}<+\infty\ ,
\ee
which is obviously weaker than (\ref{II.r7-ann}). While the ``Struwe entropy estimate'', that can be achieved for min-max problems, is ensuring 
\be
\label{II.r6-ann}
\lim_{p\rightarrow 2} \,(p-2)\,\log\delta^{-1}=0\ ,
\ee
see \cite{Lam}. The estimate (\ref{II.r.4-ann}) and the neck length control is not true in general (see a counter-exemple  for $p $-harmonic maps in \cite{LiWa-2}). In the case of $p $-harmonic maps into spheres nevertheless the situation is improving and the necklessness property holds (see \cite{LiZh} where the combination of the existence of conservation laws - i.e. (\ref{IV.9}) - with the use of Lorentz spaces, in the line of what has been introduced twenty years before in \cite{LiRi} and systematized in \cite{LR}, is central in the argument). For  general $\Om $-systems (\ref{0.6}), if one knows that the conservation law (\ref{I.5}) can be extended throughout the bubbles then theorem~\ref{th-IV.1} is implying the energy quantization and the necklessness property. However, even for $p=2$, the $(A,B)$-conservation laws are not extendable in general (see a counter-example in \cite{LR}). 

\medskip

Finally, in the case\footnote{This case covers all the $p $-harmonic approximation of coercive conformally invariant Lagrangian with minimal regularity assumptions (See \cite{Riv}). }where $|\Om|\le C_0\,|\nabla u|$  the condition~\ref{II.r.3-ann} is superfluous  thanks to the epsilon regularity theorem~\ref{th-eps-reg} assuming that we are in a neck region satisfying
\be
\label{intro-1}
\sup_{\delta<r<2^{-1}}\int_{B_{2r}(0)\setminus B_r(0)}|\nabla u |^2\ dx^2\le \ep\ ,
\ee
for $\ep>0$ small enough depending only on $m$ and $C_0$.

\medskip

To conclude, we would like to stress the potential utility of $(A,B)$ conservation laws (\ref{I.5}) for extending the main results of \cite{DGR} about the Morse index stability (question iii))
to the $p $-harmonic approximation. Such a result could be very useful for the min-max constructions of minimal or prescribed mean curvature surfaces in arbitrary closed manifolds through Sacks-Uhlenbeck relaxation.

\section{$SO(m)$ Gauge Extraction}

Let $p> 2$, $u\in W^{1,p}(\Di^2,{\R}^m)$ and $\Om\in L^2(\Di^2,so(m))$ such that
\be
\label{II.1}
\int_{\Di^2}(1+|\nabla u|^2)^{\frac{p}{2}-1}\ |\Om|^2\ dx^2<+\infty\ .
\ee
We assume that $(u,\Om)$ satisfy the following equation
\be
\label{II.2}
-\,\mbox{div}\lf((1+|\nabla u|^2)^{\frac{p}{2}-1}\,\nabla u\rg)=(1+|\nabla u|^2)^{\frac{p}{2}-1}\,\Om\cdot \nabla u\ .
\ee
Observe that under the above assumptions
\be
\label{II.3}
(1+|\nabla u|^2)^{\frac{p}{2}-1}\,\Om\, \nabla u=(1+|\nabla u|^2)^{\frac{p}{4}-\frac{1}{2}}\,\Om\cdot(1+|\nabla u|^2)^{\frac{p}{4}-\frac{1}{2}}\,\nabla u\in L^2(\Di^2)\cdot L^2(\Di^2)\hookrightarrow L^1(\Di^2)\ .
\ee
Hence we have a-priori
\be
\label{II.3-a}
\mbox{div}\lf((1+|\nabla u|^2)^{\frac{p}{2}-1}\,\nabla u\rg)\in  L^1(\Di^2)\ .
\ee
Let $Q\in W^{1,2}(\Di^2,SO(m))$ such that
\be
\label{II.4}
\int_{\Di^2}(1+|\nabla u|^2)^{\frac{p}{2}-1}\ |\nabla Q|^2\ dx^2<+\infty\ .
\ee
We have
\be
\label{II.5}
-\mbox{div}\lf((1+|\nabla u|^2)^{\frac{p}{2}-1}\,Q\,\nabla u\rg)=(1+|\nabla u|^2)^{\frac{p}{2}-1}\,\lf(Q\,\Om\,Q^{-1}-\nabla Q\,Q^{-1}\rg)\cdot Q\,\nabla u\ .
\ee
We shall be proving the following lemma which is based on the variational strategy introduced by F.H\'elein in \cite{Hel} (section 4.1) and later on adopted in a similar context by A.Schikorra in \cite{Shi}.
\begin{Lm}
\label{lm-gauge}
Let $f\in L^{p/(p-2)}(\Di^2)$ such that $f\ge 1$ and let $\Om\in L^2(\Di^2,so(m))$. Assuming 
\be
\label{II.6}
\int_{\Di^2}f\ |\Om|^2\ dx^2<+\infty\ ,
\ee
then there exists $Q\in W^{1,2}(\Di^2,SO(m))$ minimizing
\be
\label{II.7}
E_f(\Om):=\inf\lf\{  \int_{\Di^2}f\ |Q\,\Om-\nabla Q|^2\ dx^2\quad;\quad Q\in W^{1,2}_f(\Di^2,SO(m)) \quad \mbox{ and }\quad Q=Id\mbox{ on }\p \Di^2\rg\}\ ,
\ee
where
\be
\label{II.8}
W^{1,2}_f(\Di^2,SO(m)):=\inf\lf\{  Q\in W^{1,2}(\Di^2,SO(m)) \ ;\ \int_{\Di^2}f\ |\nabla Q|^2\ dx^2<+\infty\rg\}\ .
\ee
Moreover $Q$ satisfies
\be
\label{II.9}
\mbox{div}\lf(f\,\lf(Q\,\Om\,Q^{-1}-\nabla Q\,Q^{-1}\rg)  \rg)=0\ ,
\ee
and the following estimates holds
\be
\label{II.10}
\int_{\Di^2}f\ |Q\,\Om-\nabla Q|^2\ dx^2\le \int_{\Di^2}f\ |\Om|^2\ dx^2\ .
\ee
\hfill $\Box$
\end{Lm}
\noindent{Proof of lemma~\ref{lm-gauge}.} Let $Q_k$ be a minimizing sequence. We have by minimality
\be
\label{II.10-a}
\limsup_{k\rightarrow +\infty}\, \int_{\Di^2}f\ |Q_k\,\Om-\nabla Q_k|^2\ dx^2\le \int_{\Di^2}f\ |\Om|^2\ dx^2\ .
\ee
This gives in particular
\be
\label{II.11}
\begin{array}{l}
\ds\limsup_{k\rightarrow +\infty}\, \int_{\Di^2} |\nabla Q_k|^2\ dx^2\le \limsup_{k\rightarrow +\infty}\, \int_{\Di^2}f\ |\nabla Q_k|^2\ dx^2\\[5mm]
\ds\quad \le 2\,\int_{\Di^2}f\ |\Om|^2\ dx^2+2\,\limsup_{k\rightarrow +\infty}\, \int_{\Di^2}f\ |Q_k\,\Om-\nabla Q_k|^2\ dx^2<4\, \int_{\Di^2}f\ |\Om|^2\ dx^2\ .
\end{array}
\ee
We can extract a subsequence that we keep denoting $Q_k$ such that 
\be
\label{II.12}
Q_k\rightharpoonup\, Q_\infty\quad\mbox{ in }\quad W^{1,2}(\Di^2)\ ,
\ee
and by Rellich Kondrachov theorem we deduce that
\be
\label{II.13}
Q_k\rightarrow\, Q_\infty\quad\mbox{ strongly in }\quad L^q(\Di^2)\quad\mbox{ for any }q<+\infty\ .
\ee
In particular we deduce
\be
\label{II.14}
 Q_\infty \in W^{1,2}(\Di^2,SO(m))\quad\mbox{ and }\quad Q_\infty=Id_m\ .
\ee 
Let $M>1$ be arbitrary and denote, for $x\in \Di^2,$ $f_M(x):=\min\{f(x),M\}$. For any $M<+\infty$ we  obviously have   
\be
\label{II.15}
\limsup_{k\rightarrow +\infty}\, \int_{\Di^2} f_M |\nabla Q_k|^2\ dx^2\le \limsup_{k\rightarrow +\infty}\, \int_{\Di^2}f\ |\nabla Q_k|^2\ dx^2\ ,
\ee
and
\be
\label{II.16}
\limsup_{k\rightarrow +\infty}\, \int_{\Di^2}f_M\ |Q_k\,\Om-\nabla Q_k|^2\ dx^2\le E_f(\Om)\ .
\ee
We have for any $Y\in L^2(\Di^2, {R}^2\otimes M_m(R))$ and any $M\ge 1$
\be
\label{II.17}
\lim_{k\rightarrow +\infty}\int_{\Di^2}\sqrt{f_M}\, \nabla Q_k\cdot Y\ dx^2=\lim_{k\rightarrow +\infty}\int_{\Di^2} \nabla Q_k\cdot \sqrt{f_M}\,Y\ dx^2=\int_{\Di^2} \nabla Q_\infty\cdot \sqrt{f_M}\,Y\ dx^2\ .
\ee
Hence we deduce
\be
\label{II.18}
\sqrt{f_M}\, \nabla Q_k\ \rightharpoonup\ \sqrt{f_M}\, \nabla Q_\infty\ \mbox{ weakly in }L^2(\Di^2)\ ,
\ee
and thanks to Banach Steinhaus theorem and the Hilbert nature of $L^2$ we have for any $M>0$
\be
\label{II.19}
\int_{\Di^2} f_M |\nabla Q_\infty|^2\ dx^2\le \liminf_{k\rightarrow +\infty}\, \int_{\Di^2} f_M |\nabla Q_k|^2\ dx^2\le 4\ \int_{\Di^2}f\ |\Om|^2\ dx^2\ ,
\ee
and in a similar way for any $M>0$ there holds
\be
\label{II.20}
\int_{\Di^2}f_M\ |Q_\infty\,\Om-\nabla Q_\infty|^2\ dx^2\le\liminf_{k\rightarrow +\infty}\, \int_{\Di^2}f_M\ |Q_k\,\Om-\nabla Q_k|^2\ dx^2\le E_f(\Om)\ .
\ee
Since $ |Q_\infty\,\Om-\nabla Q_\infty|^2\ dx^2$ is absolutely continuous with respect to Lebesgue and since $f$ is an $L^1$ function we deduce
\be
\label{II.21}
\lim_{M\rightarrow +\infty}\int_{\Di^2}f_M\ |Q_\infty\,\Om-\nabla Q_\infty|^2\ dx^2=\int_{\Di^2}f\ |Q_\infty\,\Om-\nabla Q_\infty|^2\ dx^2
\ee
and combining (\ref{II.14}), (\ref{II.20}) together with (\ref{II.21}) we deduce
\be
\label{II.22}
\int_{\Di^2}f\ |Q_\infty\,\Om-\nabla Q_\infty|^2\ dx^2=E_f(\Om)\ ,
\ee
and $Q_\infty$ is a minimizer of (\ref{II.7}). Let $U\in C^\infty_0(\Di^2,so(m))$ and $t\in (0,1)$. Denote $Q_\infty(t):= \exp(t\,U)\, Q_\infty$. The fundamental principles of the calculus of variations imply
\be
\label{II.23}
\lf.\frac{d}{dt}\int_{\Di^2}f\ |Q_\infty(t)\,\Om-\nabla Q_\infty(t)|^2\ dx^2\rg|_{t=0}=0\ .
\ee
This gives
\be
\label{II.24}
\int_{\Di^2}\ f\, \lf<Q_\infty\,\Om-\nabla Q_\infty, U\,Q_\infty\,\Om-\nabla( U\, Q_\infty)\rg>\ dx^2=0\ .
\ee
Observe that
\be
\label{II.25}
\begin{array}{l}
\ds \lf<Q_\infty\,\Om-\nabla Q_\infty, U\,Q_\infty\,\Om-\nabla( U\, Q_\infty)\rg>=\lf<Q_\infty\,\Om Q_\infty^{-1}-\nabla Q_\infty\,Q_\infty^{-1}, U\,Q_\infty\,\Om\,Q_\infty^{-1}\rg>\\[5mm]
\ds\quad-\lf<Q_\infty\,\Om Q_\infty^{-1}-\nabla Q_\infty\,Q_\infty^{-1}, U\,\nabla Q_\infty\,Q_\infty^{-1}\rg>-\lf<Q_\infty\,\Om\,Q_\infty^{-1}-\nabla Q_\infty\,Q_\infty^{-1},\nabla U\rg>\\[5mm]
\ds\quad=\lf<Q_\infty\,\Om\, Q_\infty^{-1}-\nabla Q_\infty\,Q_\infty^{-1}, U\,\lf(Q_\infty\,\Om\,Q_\infty^{-1}-\nabla Q_\infty\,Q_\infty^{-1}\rg)\rg>\\[5mm]
\ds\quad\ -\lf<Q_\infty\,\Om\,Q_\infty^{-1}-\nabla Q_\infty\,Q_\infty^{-1},\nabla U\rg>\ .
\end{array}
\ee
Observe that for any pair of antisymmetric matrices $U$ and $V$ there holds
\be
\label{II.26}
\lf< V,U\,V\rg>=\sum_{i,k=1}^m V_{ik}\, (U\,V)_{ik}=\sum_{i,k,l=1}^m V_{ik}\, U_{il}\, V_{lk}=\frac{1}{2}\sum_{i,k,l=1}^m \lf(V_{ik}\, U_{il}\, V_{lk}+V_{lk}\, U_{li}\,V_{ik}\rg)=0\ .
\ee
Thus 
\be
\label{II.27}
\begin{array}{l}
\ds \lf<Q_\infty\,\Om-\nabla Q_\infty, U\,Q_\infty\,\Om-\nabla( U\, Q_\infty)\rg>=-\lf<Q_\infty\,\Om\,Q_\infty^{-1}-\nabla Q_\infty\,Q_\infty^{-1},\nabla U\rg>\ .
\end{array}
\ee
Hence, combining (\ref{II.24}), and (\ref{II.27}) we obtain
\be
\label{II.28}
\forall\, U\in C^\infty_0(\Di^2,so(m))\quad\quad\int_{\Di^2}\ f\, \lf<Q_\infty\,\Om\,Q_\infty^{-1}-\nabla Q_\infty\,Q_\infty^{-1},\nabla U\rg>\ dx^2=0\ .
\ee
Let now $W\in C^\infty_0(\Di^2,M_{m}({\R}))$ be arbitrary. Thanks to the  orthogonality of symmetric and antisymmetric matrices for the standard scalar product in $M_{m}({\R})$ we have
\be
\label{II.29}
\lf<Q_\infty\,\Om\,Q_\infty^{-1}-\nabla Q_\infty\,Q_\infty^{-1},\nabla W\rg>=\frac{1}{2}\lf<Q_\infty\,\Om\,Q_\infty^{-1}-\nabla Q_\infty\,Q_\infty^{-1},\nabla (W-W^t)\rg>\ .
\ee
We then deduce
\be
\label{II.30}
\forall\, W\in C^\infty_0(\Di^2,so(m))\quad\quad\int_{\Di^2}\ f\, \lf<Q_\infty\,\Om\,Q_\infty^{-1}-\nabla Q_\infty\,Q_\infty^{-1},\nabla W\rg>\ dx^2=0\ .
\ee
This identity is implying 
\be
\label{II.31}
\mbox{div}\lf(f\,\lf(Q_\infty\,\Om\,Q_\infty^{-1}-\nabla Q_\infty\,Q^{-1}_\infty\rg)  \rg)=0\ ,
\ee
and the lemma~\ref{lm-gauge} is proved.\hfill $\Box$
\par
\bigskip

Let $Q$ be the Gauge in $W^{1,2}(\Di^2,SO(m))$ given by the previous lemma for 
\be
\label{II.32}
f:=(1+|\nabla u|^2)^{\frac{p}{2}-1}\ .
\ee
We then have the existence of $\xi$ taking values into $so(m)$ such that
\be
\label{II.33}
\nabla^\perp\xi:=(1+|\nabla u|^2)^{\frac{p}{2}-1}\,\lf(Q\,\Om\,Q^{-1}-\nabla Q\,Q^{-1}\rg)\ .
\ee
Observe that we have respectively
\be
\label{II.34}
\begin{array}{l}
\ds\int_{\Di^2}|\nabla \xi|^{\frac{p}{p-1}}\ dx^2=\int_{\Di^2} \lf|(1+|\nabla u|^2)^{\frac{p}{4}-\frac{1}{2}}\,\lf(Q\,\Om\,Q^{-1}-\nabla Q\,Q^{-1}\rg)\rg|^{\frac{p}{p-1}}\,\lf|(1+|\nabla u|^2)^{\frac{p}{4}-\frac{1}{2}}\rg|^{\frac{p}{p-1}}\, dx^2\\[5mm]
\ds\quad \le \lf[ \int_{\Di^2}(1+|\nabla u|^2)^{\frac{p}{2}-1}\,|Q\,\Om-\nabla Q|^2   \rg]^{\frac{p}{2(p-1)}}\ \lf[\int_{\Di^2}\lf|(1+|\nabla u|^2)^{\frac{p}{4}-\frac{1}{2}}\rg|^{\frac{p}{p-1}\frac{2p-2}{p-2}}\, dx^2\rg]^{\frac{p-2}{2p-2}}\\[5mm]
\ds\quad\le \lf[ \int_{\Di^2}(1+|\nabla u|^2)^{\frac{p}{2}-1}\,|Q\,\Om-\nabla Q|^2   \rg]^{\frac{p}{2(p-1)}}\ \lf[\int_{\Di^2}(1+|\nabla u|^2)^{\frac{p}{2}}\, dx^2\rg]^{\frac{p-2}{2p-2}}<+\infty\ .
\end{array}
\ee
and
\be
\label{II.35}
\int_{\Di^2}\frac{|\nabla \xi|^2}{(1+|\nabla u|^2)^{\frac{p}{2}-1}}\ dx^2\le  \int_{\Di^2}(1+|\nabla u|^2)^{\frac{p}{2}-1}\,|Q\,\Om-\nabla Q|^2\  dx^2<+\infty\ .
\ee
With this notation the system (\ref{II.5}) becomes
\be
\label{II.36}
-\mbox{div}\lf((1+|\nabla u|^2)^{\frac{p}{2}-1}\,Q\,\nabla u\rg)=\nabla^\perp\xi\cdot Q\,\nabla u\ .
\ee
In the next section we are going to ``refine'' the gauge change $Q$ in the spirit of what is made in \cite{Riv} in order to obtain a conservation law.
\section {The $Gl_m$ Gauge and the $(A,B)$ Conservation Law.}
\reset
{\bf Proof of theorem~\ref{th-I.1}.} Let $\Om\in L^2(\Di^2,{\R}^2\otimes so(m))$ and $u\in W^{1,p}(\Di^2,{\R}^m)$ satisfying
 \be
\label{III.1}
\int_{\Di^2}(1+|\nabla u|^2)^{\frac{p}{2}-1}\ |\Om|^2\ dx^2<\sigma\ ,
\ee
for some $\sigma>0$ to be fixed later on. We assume that $(u,\Om)$ satisfy the following equation
\be
\label{III.2}
-\,\mbox{div}\lf((1+|\nabla u|^2)^{\frac{p}{2}-1}\,\nabla u\rg)=(1+|\nabla u|^2)^{\frac{p}{2}-1}\,\Om\cdot \nabla u\ .
\ee
Thanks to the previous section we have the existence of $Q\in W^{1,2}(\Di^2,SO(m))$ such that
\be
\label{III.3}
\mbox{div}\lf((1+|\nabla u|^2)^{\frac{p}{2}-1}\,\lf(Q\,\Om\,Q^{-1}-\nabla Q\,Q^{-1}\rg)  \rg)=0\ ,
\ee
and the following estimates holds
\be
\label{III.4}
\int_{\Di^2}(1+|\nabla u|^2)^{\frac{p}{2}-1} |Q\,\Om-\nabla Q|^2\ dx^2\le \int_{\Di^2}(1+|\nabla u|^2)^{\frac{p}{2}-1} |\Om|^2\ dx^2<\sigma\ .
\ee
and introducing
\be
\label{III.5}
\nabla^\perp\xi:=(1+|\nabla u|^2)^{\frac{p}{2}-1}\,\lf(Q\,\Om\,Q^{-1}-\nabla Q\,Q^{-1}\rg)\ ,
\ee
there holds
\be
\label{III.6}
-\mbox{div}\lf((1+|\nabla u|^2)^{\frac{p}{2}-1}\,Q\,\nabla u\rg)=\nabla^\perp\xi\cdot Q\,\nabla u\ .
\ee
We are looking for $\ep\in L^\infty\cap W^{1,2}(\Di^2,M_m({\R}))$ such that
\be
\label{III.7}
\int_{\Di^2}(1+|\nabla u|^2)^{\frac{p}{2}-1}\ |\nabla\ep|^2\ dx^2<C\,\int_{\Di^2}(1+|\nabla u|^2)^{\frac{p}{2}-1}\ |\Om|^2\ dx^2\ ,
\ee
where $C>0$ is universal and 
\be
\label{III.8}
\mbox{div}\lf( (Id_m+\ep)\,\nabla^\perp\xi\,Q-(1+|\nabla u|^2)^{\frac{p}{2}-1}\, \nabla\,\ep\,Q\,\rg)=0\ .
\ee
Introducing
\be
\label{III.9}
-\nabla^\perp B:=-(1+|\nabla u|^2)^{\frac{p}{2}-1}\, \nabla\,\ep\,Q+ (Id_m+\ep)\,\nabla^\perp\xi\,Q\ .
\ee
In particular this implies 
\be
\label{III.10}
\begin{array}{rl}
\ds-\mbox{div}\lf((1+|\nabla u|^2)^{\frac{p}{2}-1}\,(Id_m+\ep)\,Q\,\nabla u\rg)&\ds=\lf[-(1+|\nabla u|^2)^{\frac{p}{2}-1}\, \nabla\,\ep+ (Id_m+\ep)\,\nabla^\perp\xi\rg]\,Q\,\nabla u\\[5mm]
 &\ds=\mbox{div}(B\cdot\nabla^\perp u)\ .
\end{array}
\ee
\noindent{\bf Construction of $\ep$.}  Denote
\be
\label{III.11}
f(x):=(1+|\nabla u|^2(x))^{\frac{p}{2}-1}\quad\mbox{ and }\quad f_M(x):=\min\{f(x),M\}\ .
\ee
We introduce $Q_M$ minimizing $E_{f_M}(\Om)$ given by lemma~\ref{lm-gauge} and $\xi_M$ satisfying
\be
\label{III.11-a}
\nabla^\perp\xi_M:=f_M\,\lf(Q_M\,\Om\,Q_M^{-1}-\nabla Q_M\,Q_M^{-1}\rg)\ ,
\ee
Because of (\ref{II.35}) we have
and
\be
\label{III.11-b}
\begin{array}{l}
\ds\int_{\Di^2}\frac{|\nabla \xi_M|^2}{f_M}\ dx^2\le  \int_{\Di^2}f_M\,|Q_M\,\Om-\nabla Q_M|^2\  dx^2\\[5mm]
\ds\quad\le \int_{\Di^2}f_M\,|\Om|^2\  dx^2\le \int_{\Di^2}(1+|\nabla u|^2)^{\frac{p}{2}-1}\ |\Om|^2\ dx^2<\sigma\ ,
\end{array}
\ee
together with
\be
\label{III.11-c}
\int_{\Di^2}f_M\,|\nabla Q_M|^2\  dx^2\le 2\,  \int_{\Di^2}f_M\,[|Q_M\,\Om-\nabla Q_M|^2+|\Om|^2]\  dx^2\le 2\, \sigma\ .
\ee
Moreover, arguing as in (\ref{II.34}) we have
\be
\label{III.11-d}
\int_{\Di^2}|\nabla\xi_M|^{p/(p-1)}\ dx^2\le \lf[ \int_{\Di^2}f_M\,|Q_M\,\Om-\nabla Q_M|^2\  dx^2\rg]^{\frac{2p}{2(p-1)}}\ \lf[\int_{\Di^2}(1+|\nabla u|^2)^{\frac{p}{2}}\ dx^2\rg]^{\frac{p-2}{2p-2}}\ .
\ee
We are then looking for $(\e_M,B_M)$  solving
\be
\label{III.12}
\lf\{
\begin{array}{l}
\ds-\mbox{div}\lf(f_M \,\nabla\,\ep_M\rg)=-\mbox{ div}\lf(\ep_M\,\nabla^\perp\xi_M\rg)+\nabla B_M\cdot\nabla^\perp Q_M^{-1}  \quad\mbox{ in }\Di^2\\[5mm]
\ds\ep_M=0\quad\mbox{ on }\p \Di^2\ ,
\end{array}
\rg.
\ee
and
\be
\label{III.12-a}
\lf\{
\begin{array}{l}
\ds\mbox{ div}\lf(\frac{\nabla B_M}{{f_M}} \rg)= -\mbox{ div}(\nabla^\perp \ep_M\,Q_M)- \mbox{ div}\lf( (Id_m+\ep_M) \frac{\nabla\xi_M}{f_M}\,Q_M \rg) \\[5mm]
\ds\p_\nu B_M:=-\,\p_\nu\xi\quad\mbox{ on }\p \Di^2
\end{array}
\rg.
\ee
Observe that, $\ep_M$ being given, $B_M$ is the minimizer in $W^{1,2}(\Di^2,M_m({\R}^m))$ of 
\be
\label{III.13}
\int_{\Di^2} \lf|\frac{\nabla B}{\sqrt{f}_M}+\sqrt{f_M}\nabla^\perp\ep_M\, Q_M+(Id_m+\ep_M) \frac{\nabla\xi_M}{\sqrt{f_M}}\,Q_M \rg|^2\ dx^2\ .
\ee
Indeed, the minimizer satisfies
\be
\label{III.14}
\begin{array}{l}
\ds\forall\, \beta\in W^{1,2}(\Di^2,M_m({\R}^m))\\[5mm]
\ds\int_{\Di^2}\frac{\nabla\beta}{\sqrt{f_M}}\cdot \lf[\frac{\nabla B_M}{\sqrt{f}_M}+\sqrt{f_M}\nabla^\perp\ep_M\, Q_M+(Id_m+\ep_M) \frac{\nabla\xi_M}{\sqrt{f_M}}\,Q_M \rg]\ dx^2=0\ .
\end{array}
\ee
While $B_M$ being given, $\ep_M$ solving (\ref{III.12}) is given by lemma~\ref{lm-wente-weight}. We then construct $(\ep_M,B_M)$ simultaneously by iteration.\par We define
\be
\label{III.15}
\lf\{
\begin{array}{l}
\ds\mbox{ div}\lf(\frac{\nabla B_{M,0}}{{f_M}} \rg)=- \mbox{ div}\lf( \frac{\nabla\xi_M}{f_M}\,Q_M \rg) \\[5mm]
\ds\p_\nu B_{M,0}:=-\,\p_\nu\xi_M\quad\mbox{ on }\p \Di^2
\end{array}
\rg.
\ee
and $\ep_{M,0}=0$. Multipying by $B_{M,0}$ and integrating by parts is giving
\be
\label{III.15-a}
-\,\int_{\Di^2}\frac{|\nabla B_{M,0}|^2}{{f_M}}\ dx^2=\int_{\Di^2}f_M^{-1}\lf<\nabla B_{M,0}\cdot \nabla\xi_M\,Q_M\rg>\ dx^2\ .
\ee
Using Cauchy Schwartz this gives
\be
\label{III.15-b}
\int_{\Di^2}\frac{|\nabla B_{M,0}|^2}{{f_M}}\ dx^2\le\int_{\Di^2}\frac{|\nabla\xi_M|^2}{f_M}\ dx^2\ .
\ee
Then we define iteratively for $k\ge 1$
\be
\label{III.16}
\lf\{
\begin{array}{l}
\ds-\mbox{div}\lf(f_M \,\nabla\,\ep_{M,k}\rg)=-\mbox{ div}\lf(\ep_{M,k-1}\,\nabla^\perp\xi_M\rg)+\nabla B_{M,k-1}\cdot\nabla^\perp Q^{-1}_M  \quad\mbox{ in }\Di^2\\[5mm]
\ds\ep_{M,k}=0\quad\mbox{ on }\p \Di^2\ ,
\end{array}
\rg.
\ee
and
\be
\label{III.17}
\lf\{
\begin{array}{l}
\ds\mbox{ div}\lf(\frac{\nabla B_{M,k}}{{f_M}} \rg)=- \mbox{ div}(\nabla^\perp \ep_{M,k-1}\,Q_M)- \mbox{ div}\lf( \ep_{M,k-1} \frac{\nabla\xi_M}{f_M}\,Q_M \rg) \\[5mm]
\ds\p_\nu B_{M,k}=0\quad\mbox{ on }\p \Di^2
\end{array}
\rg.
\ee
Because of lemma~\ref{lm-wente-weight} we have for $k\ge 0$ the estimate 
\be
\label{III.18}
\begin{array}{l}
\ds\|\ep_{M,k}\|_\infty^2+\int_{\Di^2} f_M\, |\nabla \ep_{M,k}|^2\ dx^2\le C_0\ \int_{\Di^2}\frac{|\nabla\xi_M|^2}{f_M}\ dx^2\   \int_{\Di^2} f_M\, |\nabla \ep_{M,k-1}|^2\ dx^2\\[5mm]
\ds\quad+C_0\   \int_{\Di^2} f_M\, |\nabla Q_M|^2\ dx^2\ \int_{\Di^2}\frac{|\nabla B_{M,k-1}|^2}{f_M}\ dx^2\ .
\end{array}
\ee
Multiplying (\ref{III.17}) by $B_{M,k}$ and integrating is giving
\be
\label{III.19}
\begin{array}{l}
\ds\int_{\Di^2}\frac{|\nabla B_{M,k}|^2}{f_M}\ dx^2\le \|\ep_{M,k-1}\|_\infty\ \sqrt{\int_{\Di^2} f_M\, |\nabla Q_M|^2\ dx^2\ \int_{\Di^2}\frac{|\nabla B_{M,k}|^2}{f_M}\ dx^2}\\[5mm]  
\ds +\|\ep_{M,k-1}\|_\infty\ \sqrt{\int_{\Di^2} \frac{ |\nabla \xi_M|^2}{f_M}\ dx^2\ \int_{\Di^2}\frac{|\nabla B_{M,k}|^2}{f_M}\ dx^2}\ .
\end{array}
\ee
We deduce
\be
\label{III.20}
\begin{array}{l}
\ds\int_{\Di^2}\frac{|\nabla B_{M,k}|^2}{f_M}\ dx^2\le 2\ \|\ep_{M,k-1}\|^2_\infty\ \lf[\int_{\Di^2} f_M\, |\nabla Q_M|^2+ \frac{ |\nabla \xi_M|^2}{f_M} dx^2\rg]\ .
\end{array}
\ee
Combining (\ref{III.11-b}), (\ref{III.11-c}), (\ref{III.18}) and (\ref{III.20}) gives then the existence of a universal constant $C_1>0$ such that $\forall k\ge 1$
\be
\label{III.22}
\begin{array}{l}
\ds\|\ep_{M,k}\|_\infty^2+\int_{\Di^2} f_M\, |\nabla \ep_{M,k}|^2+\frac{|\nabla B_{M,k}|^2}{f_M}\ dx^2\\[5mm]
\ds\quad\le C_1\, \sigma\ \|\ep_{M,k-1}\|_\infty^2+ C_1\, \sigma\ \int_{\Di^2} f_M\, |\nabla \ep_{M,k-1}|^2+\frac{|\nabla B_{M,k-1}|^2}{f_M}\ dx^2\ .
\end{array}
\ee
Iterating over $k$ gives
\be
\label{III.23}
\begin{array}{l}
\ds\|\ep_{M,k}\|_\infty+\sqrt{\int_{\Di^2} f_M\, |\nabla \ep_{M,k}|^2}+\sqrt{\int_{\Di^2}\,\frac{|\nabla B_{M,k}|^2}{f_M}\ dx^2}\\[5mm]
\ds\quad\le (\sqrt{C_1\,\sigma})^k\ \sqrt{\int_{\Di^2}\frac{|\nabla B_{M,0}|^2}{{f_M}}\ dx^2}\le  (\sqrt{C_1\,\sigma})^k\ \sqrt{\int_{\Di^2}\frac{|\nabla\xi_M|^2}{f_M}\ dx^2}\ ,
\end{array}
\ee
that we can also rewrite as
\be
\label{III.24}
\|\ep_{M,k}\|_\infty+\|\sqrt{f_M}\nabla\ep_{M,k}\|_{L^2(\Di^2)}+\|\sqrt{f_M}^{-1}\nabla B_{M,k}\|_{L^2(\Di^2)}\le  (\sqrt{C_1\,\sigma})^k\ \|\sqrt{f_M}^{-1}\nabla\xi_M \|_{L^2(\Di^2)}\ .
\ee
Choosing $C_1\sigma<1$ arbitrary we have using the triangular inequality
\be
\label{III.25}
\begin{array}{l}
\ds\lf\|\sum_{k=0}^{+\infty}\ep_{M,k}\rg\|_\infty+\lf\|\sqrt{f_M}\sum_{k=0}^{+\infty}\nabla\ep_{M,k}\rg\|_{L^2(\Di^2)}+\lf\|\sqrt{f_M}^{-1} \sum_{k=0}^{+\infty}\nabla B_{M,k}\rg\|_{L^2(\Di^2)}\\[5mm]
\ds\quad\le C_2\ \|\sqrt{f_M}^{-1}\nabla\xi_M \|_{L^2(\Di^2)}\ .
\end{array}
\ee
This bound permits to sum up the equations (\ref{III.16}) and (\ref{III.17}) in order to obtain that 
\[
(\ep_M,B_M):=\lf(\sum_{k=0}^{+\infty}\ep_{M,k},\sum_{k=0}^{+\infty} B_{M,k}\rg)
\]
is a solution of the system (\ref{III.12}) and (\ref{III.12-a}) such that
\be
\label{III.26}
\begin{array}{l}
\ds\lf\|\ep_{M}\rg\|_\infty+\lf\|\sqrt{f_M}\,\nabla\ep_{M}\rg\|_{L^2(\Di^2)}+\lf\|\sqrt{f_M}^{-1}\,\nabla B_{M}\rg\|_{L^2(\Di^2)}\le C_2\ \|\sqrt{f_M}^{-1}\nabla\xi_M \|_{L^2(\Di^2)}\ .
\end{array}
\ee
Let $D_M$ satisfy
\be
\label{III.27}
\nabla^\perp B_M-f_M\,\nabla\ep_M\, Q_M+(Id_m+\ep_M)\nabla^\perp\xi_M\,Q_M=f_M\,\nabla D_M\ .
\ee
Since $\ep_M=0$ on $\p \Di^2$, $Q_M=Id_m$ on $\p \Di^2$ and $\p_\nu B_M+\p_\nu\xi_M=0$ on $\p \Di^2$, we have $\p_\tau D_M=0$ on $\p \Di^2$ (where $\tau$ is the unit tangent to the boundary of $\Di^2$) and we can choose $D_M=0$ on $\p \Di^2$. Multiplying (\ref{III.27}) on the left by $Q_M^{-1}$ and using (\ref{III.12}) gives
\be
\label{III.28}
\lf\{
\begin{array}{l}
\ds\mbox{div}\lf(f_M\,\nabla D_M\,Q_M^{-1}   \rg)=0\quad\mbox{ in }\Di^2\\[5mm]
\ds D_M=0\quad\mbox{ on }\p \Di^2\ .
\end{array}
\rg.
\ee
We perform the following Hodge decomposition of $\nabla D_M\,Q_M^{-1}$ :
\be
\label{III.29}
\nabla D_M\,Q_M^{-1} =\nabla \zeta_M+\frac{\nabla^\perp\eta_M}{f_M}\ ,
\ee
where $ \zeta_M$ is the minimizer of
\be
\label{III.30}
\inf\lf\{ \int_{\Di^2}f_M\, |\nabla\zeta- \nabla D_M\,Q_M^{-1}|^2\ dx^2\ ;\ \zeta\in W^{1,2}_0(\Di^2,M_m({\R}))\rg\}\ .
\ee
Because of (\ref{III.28}) $\zeta_M$ satisfies
\be
\label{III.31}
\lf\{
\begin{array}{l}
\ds\mbox{div}\lf(f_M\,\nabla \zeta_M   \rg)=0\quad\mbox{ in }\Di^2\\[5mm]
\ds \zeta_M=0\quad\mbox{ on }\p \Di^2\ .
\end{array}
\rg.
\ee
Multiplying by $\zeta_M$ and integrating by parts is giving $\zeta_M\equiv 0$. Hence
\be
\label{III.32}
\nabla D_M\,Q_M^{-1} =\frac{\nabla^\perp\eta_M}{f_M}\ .
\ee
This gives that $D_M$ satisfies
\be
\label{III.33}
\lf\{
\begin{array}{l}
\ds\mbox{div}\lf(f_M\,\nabla D_M   \rg)=\nabla^\perp\eta_M\cdot \nabla Q_M \quad\mbox{ in }\Di^2\\[5mm]
\ds D_M=0\quad\mbox{ on }\p \Di^2\ .
\end{array}
\rg.
\ee
Using again  lemma~\ref{lm-wente-weight} we obtain the following inequality
\be
\label{III.34}
\begin{array}{l}
\ds\int_{\Di^2}f_M\,|\nabla D_M|^2\ dx^2\le C_0\, \int_{\Di^2}f_M\,|\nabla Q_M|^2\ dx^2\ \int_{\Di^2}\frac{|\nabla\eta_M|^2}{f_M}\ dx^2\\[5mm]
\ds\quad=C_0\, \int_{\Di^2}f_M\,|\nabla Q_M|^2\ dx^2\ \int_{\Di^2}f_M\,|\nabla D_M|^2\ dx^2\\[5mm]
\ds\quad\le 2\, \sigma\, C_0 \int_{\Di^2}f_M\,|\nabla D_M|^2\ dx^2\, .
\end{array}
\ee
We choose $2\, \sigma\, C_0<1$ and we obtain $D_M=0$. Hence we have finally
\be
\label{III.35}
\nabla^\perp B_M-f_M\,\nabla\ep_M\, Q_M+(Id_m+\ep_M)\nabla^\perp\xi_M\,Q_M=0\ .
\ee
We have using dominated convergence that that $\sqrt{f_M}$ is converging strongly towards $\sqrt{f}$ in $L^{2p'}(\Di^2)$ moreover $\sqrt{f_M}\nabla\ep_M$ is uniformly bounded
in $L^2$. Hence since $f_M\ge 1$ we can extract a subsequence $M_k\rightarrow+\infty$ such that we have simultaneously
\be
\label{III.36}
\sqrt{f_{M_k}}\,\nabla\ep_{M_k}\rightharpoonup X\quad\mbox{ weakly in }L^2(\Di^2)
\ee
and
\be
\label{III.37}
\nabla\ep_{M_k}\rightharpoonup \nabla\ep\quad\mbox{ weakly in }L^2(\Di^2)\ .
\ee
Since $\sqrt{f_M}$ is converging strongly towards $\sqrt{f}$ in $L^{2p'}(\Di^2)$ wee deduce that
\be
\label{III.38}
\sqrt{f_{M_k}}\,\nabla\ep_{M_k}\rightharpoonup \sqrt{f}\,\nabla\ep\quad\mbox{ weakly in }{\mathcal D}'(\Di^2)\ .
\ee
This gives $X=\sqrt{f}\,\nabla\ep$ and then
\be
\label{III.38-a}
\sqrt{f_{M_k}}\,\nabla\ep_{M_k}\rightharpoonup \sqrt{f}\,\nabla\ep\quad\mbox{ weakly in }L^2(\Di^2)\ .
\ee
We have also that, after extracting a subsequence, $Q_{M_k}$ converges weakly to $Q$ in $W^{1,2}_{Id}(\Di^2,SO(m))$. Since $\ep_{M_k}$ is also weakly converging towards $\ep$
in $W^{1,2}_{0}(\Di^2,M_m({\R})$ and since it is uniformly bounded in $L^\infty$ (thanks to (\ref{III.26}) we deduce using Rellich-Kondrachov
\be
\label{III.39}
\forall i,j,r,s\in \{1,\ldots, m\} \quad\quad \ep^{ij}_{M_k}\, Q_{M_k}^{rs}\longrightarrow \ \ep^{ij}\, Q^{rs}\quad\mbox{ strongly in }L^q\quad\forall q<+\infty\, .
\ee
We have also that
\be
\label{III.40}
f_{M_k}\,\lf(Q_{M_k}\,\Om\,Q_{M_k}^{-1}-\nabla Q_{M_k}\,Q_{M_k}^{-1}\rg)\ \rightharpoonup \ f\,\lf(Q\,\Om\,Q^{-1}-\nabla Q\,Q^{-1}\rg)\ \mbox{ weakly in }{\mathcal D}'(\Di^2)\ .
\ee
Because of (\ref{III.11-d}) we have that
\be
\label{III.41}
\nabla^\perp\xi_{M_k}\ \rightharpoonup\  \ f\,\lf(Q\,\Om\,Q^{-1}-\nabla Q\,Q^{-1}\rg)\ \mbox{ weakly in }L^{p/(p-1)}(\Di^2)\ .
\ee
We deduce that $\nabla^\perp B_{M_k}$ is uniformly bounded in $L^{2p/(2p-1)}(\Di^2)$ and it converges weakly in this space towards $B$ satisfying
\be
\label{III.42}
\nabla^\perp B-f\,\nabla\ep\, Q+(Id_m+\ep)\nabla^\perp\xi\,Q=0\
\ee
Moreover we have
\be
\label{III.43}
\begin{array}{l}
\ds\lf\|\ep\rg\|_{L^\infty(\Di^2)}+\lf\|\sqrt{f}\,\nabla\ep\rg\|_{L^2(\Di^2)}+\lf\|\sqrt{f^{-1}}\,\nabla B\rg\|_{L^2(\Di^2)}\le C_2\ \|\sqrt{f^{-1}}\nabla\xi \|_{L^2(\Di^2)}\\[5mm]
\ds \quad \le C_3\ \lf[\int_{\Di^2} (1+|\nabla u|^2)^{\frac{p}{2}-1}\ |\Om|^2\ dx^2\rg]^{1/2}\ .
\end{array}
\ee
\section{Uniform $L^{2,1}$ Estimates}
\reset
\subsection{The Disc Case}
The following result can also be obtained by establishing an epsilon regularity theorem independent of $p\ge 2$. We present a direct method that
we are going to follow for the degenerating annulus case in the next subsection, a situation where the epsilon regularity is a priori inoperative. 
\begin{Lm}
\label{th-IV.1}
There exists $\sigma>0$ such that for any $p\ge 2$ the following holds. Let $\Om\in L^2(\Di^2,{\R}^2\otimes so(m))$ and $u\in W^{1,p}(\Di^2,{\R}^m)$ satisfying
\be
\label{I-r.1}
-\,\mbox{div}\lf((1+|\nabla u|^2)^{\frac{p}{2}-1}\,\nabla u\rg)=(1+|\nabla u|^2)^{\frac{p}{2}-1}\,\Om\cdot \nabla u\ .
\ee
Assume that 
 \be
\label{I-r.2}
\int_{\Di^2}(1+|\nabla u|^2)^{\frac{p}{2}-1}\ |\Om|^2\ dx^2<\sigma\ ,
\ee
then for any $0<t<1$ there exists $C(t)>0$ depending only on $t$ such that
\be
\label{I-r-b2}
\begin{array}{l}
\|\nabla u\|_{L^{2,1}(B_t(0))}\le \  C \ (p-2)^\al\ \lf[ \|\nabla u\|_{L^p(\Di^2)}^{3/4}+ \|\nabla u\|_{L^p(\Di^2)}\rg]^{\frac{1}{p-1}}\  \lf(1+\|\nabla u\|_{L^p(\Di^2)}^2\rg)^{\frac{p-2}{2(p-1)}} \\[5mm]
\ds\ +\ C(t)\ \lf[\int_{\Di^2}\ (1+|\nabla u|^2)^{\frac{p}{2}-1}\ |\Om|^2\  dx^2\rg]^{\frac{1}{2(p-1)}}\ \lf(1+\|\nabla u\|_{L^p(\Di^2)}^2\rg)^{\frac{p}{4(p-1)}} \\[5mm]
\ds\ +\,C(t)\ \|\nabla u\|_{L^{p}(\Di^2)}^{\frac{1}{p-1}}\ \lf(1+\|\nabla u\|_{L^p(\Di^2)}^2\rg)^{\frac{p-2}{2(p-1)}}\ . 
\end{array}
\ee
\hfill $\Box$
\end{Lm}
\noindent{\bf Proof of lemma~\ref{th-IV.1}} From theorem~\ref{th-I.1} wee have the existence of $A\in L^\infty\cap W^{1,2}(\Di^2,Gl_m({\R}))$ and $B\in W^{1,p/(p-1)}(\Di^2,M_m({\R})$ such that
\be
\label{I-r.3}
\lf\|\mbox{dist}(A,SO(m))\rg\|^2_{L^\infty(\Di^2)}+\int_{\Di^2}(1+|\nabla u|^2)^{\frac{p}{2}-1}\,|\nabla A|^2\ dx^2\le C\, \int_{\Di^2}(1+|\nabla u|^2)^{\frac{p}{2}-1}\ |\Om|^2\ dx^2
\ee
such that
\be
\label{I-r.4}
\int_{\Di^2}\frac{|\nabla B|^2}{(1+|\nabla u|^2)^{\frac{p}{2}-1}\,}\ dx^2\le C\, \int_{\Di^2}(1+|\nabla u|^2)^{\frac{p}{2}-1}\ |\Om|^2\ dx^2
\ee
and equation (\ref{I-r.1}) can be rewritten as follows
\be
\label{IV.9-b}
-\,\mbox{div}\lf((1+|\nabla u|^2)^{\frac{p}{2}-1}\,A\,\nabla u\rg)=\nabla^\perp B\cdot\nabla u\quad\mbox{ in }{\mathcal D}'(\Di^2)\ .
\ee
Observe that
\be
\label{IV.12}
\begin{array}{l}
\ds\int_{\Di^2}|\nabla B|^{p/(p-1)}\ dx^2=\int_{\Di^2}\lf[\frac{|\nabla B|^2}{(1+|\nabla u|^2)^{\frac{p}{2}-1}\,}    \rg]^{\frac{p}{2(p-1)}}\ \lf[ (1+|\nabla u|^2)^{\frac{p}{2}-1}  \rg]^{\frac{p}{2(p-1)}}\ dx^2\\[5mm]
\ds\quad\quad\le \lf[ \int_{\Di^2}\frac{|\nabla B|^2}{(1+|\nabla u|^2)^{\frac{p}{2}-1}\,}\ dx^2  \rg]^{\frac{p}{2(p-1)}}\ \lf[  \int_{\Di^2}\  \lf[ (1+|\nabla u|^2)^{\frac{p}{2}-1}  \rg]^{\frac{p}{p-2}}\ dx^2 \rg]^{\frac{p-2}{2(p-1)}}\\[5mm]
\ds\quad\quad\le  \lf[ \int_{\Di^2}\ (1+|\nabla u|^2)^{\frac{p}{2}-1}\ |\Om|^2\  dx^2\rg]^{\frac{p}{2(p-1)}}\ \lf[  \int_{\Di^2}\  \lf(\sqrt{1+|\nabla u|^2}\rg)^p  dx^2 \rg]^{\frac{p-2}{2(p-1)}}\ ,
\end{array}
\ee
which gives
\be
\label{IV.12-a}
\|\nabla B\|_{L^{p'}(\Di^2)}\le C_p \ \sqrt{\int_{\Di^2}\ (1+|\nabla u|^2)^{\frac{p}{2}-1}\ |\Om|^2\  dx^2}\ \ \lf[  \int_{\Di^2}\  \lf(1+|\nabla u|^2\rg)^{\frac{p}{2}}  dx^2 \rg]^{\frac{p-2}{2p}}\ .
\ee
Let
\be
\label{IV.10}
(1+|\nabla u|^2)^{\frac{p}{2}-1}\,A\,\nabla u=A\,\nabla\varphi+\nabla^\perp\psi
\ee
where
\be
\label{IV.11}
\lf\{
\begin{array}{rl}
\ds -\mbox{ div}(A\,\nabla\varphi)&\ds=\nabla^\perp B\cdot\nabla u\quad\mbox{ in }\Di^2\\[5mm]
\ds\varphi&\ds=0\quad\mbox{ on }\p \Di^2\ .
\end{array}
\rg.
\ee
The existence and uniqueness from $\varphi$ is given by lemma~\ref{lm-wente-type-estimates} moreover
\be
\label{IV.17}
\|\nabla\varphi\|_{L^{2,1}(\Di^2)}\le C\ \|\nabla u\|_{L^{p}(\Di^2)}\  \|\nabla B\|_{L^{p'}(\Di^2)}\ ,
\ee
where $C>0$ is independent of $p\in[2,3]$. Hence (\ref{IV.12-a}) is giving
\be
\label{IV.17-a}
\|\nabla\varphi\|_{L^{2,1}(\Di^2)}\le C\, \sqrt{\int_{\Di^2}\ (1+|\nabla u|^2)^{\frac{p}{2}-1}\ |\Om|^2\  dx^2}\ \ \sqrt{ \int_{\Di^2}\  \lf(1+|\nabla u|^2\rg)^{\frac{p}{2}}  dx^2} \ ,
\ee
where $C>0$ is independent of $p\in[2,3]$.

\medskip

We introduce now two operators : For any $X\in L^p(\Di^2,{\R}^2\otimes {\R}^m)$ we define respectively
\be
\label{IV.18}
\begin{array}{rcl}
\ds S_A\ :\ L^p(\Di^2,{\R}^2\otimes {\R}^m)\ &\ds\longrightarrow\ &\ds L^{p'}(\Di^2,{\R}^2\otimes {\R}^m)\\[5mm]
    X&\ds\longrightarrow & \ds S_A(X):= \frac{(1+|A^{-1}\,X|^2)^{p/2-1}}{(1+\|A^{-1}\,X\|_{L^p(\Di^2)}^2)^{p/2-1}}\ X\ ,
    \end{array}
\ee
and for $s\in [3/2,3]$
\be
\label{IV.19}
\begin{array}{rcl}
\ds T_A\ :\  L^s(\Di^2,{\R}^2\otimes {\R}^m)\ &\ds\longrightarrow\ &\ds L^{s}(\Di^2,{\R}^2\otimes {\R}^m)\\[5mm]
 \ds   X&\ds\longrightarrow & \ds \nabla^\perp w\ ,
\end{array}
\ee
where $w$ is the unique solution in $W^{1,s}_0(\Di^2,{\R}^m)$ of
\be
\label{IV.20}
\lf\{
\begin{array}{l}
\ds -\,\mbox{div}\lf(A^{-1} \nabla w\rg)= \mbox{ div }\lf( A^{-1}\,X^\perp  \rg)\quad\mbox{ in }\Di^2\\[3mm]
\ds w=0\quad\mbox{ on }\p \Di^2\ .
\end{array}
\rg.
\ee
Observe that in particular $T_A(A\,\nabla u)=0$. We consider the commutator $S_A\circ T_A-T_A\circ S_A$ on $L^p(\Di^2,{\R}^2\otimes {\R}^m)$ and using a slight adjustment of lemma A.5 of \cite{BR} where $a(x)\,|X|$ is replaced by $|A^{-1} X|$, we obtain
\be
\label{IV.21}
\lf\| T_A\circ S_A(A\,\nabla u) \rg\|_{L^{p/(p-1)}(\Di^2)}\le C\ (p-2)\ \lf[ \|\nabla u\|_{L^p(\Di^2)}^{3/4}+ \|\nabla u\|_{L^p(\Di^2)}\rg]\ .
\ee
We also have  
\be
\label{IV.21-a}
S_A(A\,\nabla u)=\frac{\lf(1+|\nabla u|^2\rg)^{p/2-1}}{\ds\lf(1+\|\nabla u\|_{L^p(\Di^2)}^2\rg)^{p/2-1}}\ A\,\nabla u\ .
\ee
Let $w$ such that $T_A\circ S_A(A\,\nabla u)=\nabla^\perp w$. It solves
\be
\label{IV.21-b}
\lf\{
\begin{array}{l}
\ds -\,\mbox{div}\lf(A^{-1} \nabla w\rg)= \mbox{ div }\lf(  \frac{\lf(1+|\nabla u|^2\rg)^{p/2-1}}{\ds\lf(1+\|\nabla u\|_{L^p(\Di^2)}^2\rg)^{p/2-1}}\ \nabla^\perp u\rg)\quad\mbox{ in }\Di^2\\[3mm]
\ds w=0\quad\mbox{ on }\p \Di^2\ .
\end{array}
\rg.
\ee
This gives the existence of $v$ such that
\be
\label{IV.22}
\frac{\ds \lf(1+|\nabla u|^2\rg)^{p/2-1}}{\ds\lf(1+\|\nabla u\|_{L^p(\Di^2)}^2\rg)^{p/2-1}}\ \nabla u=\nabla v+A^{-1}\,\nabla^\perp w\ .
\ee
Classical elliptic estimates applied to (\ref{IV.20}) give
\be
\label{IV.23-a}
\begin{array}{rl}
\ds\lf\|\lf(1+\|\nabla u\|_{L^p(\Di^2)}^2\rg)^{p/2-1}\,\nabla w\rg\|_{L^{\frac{p}{p-1}}(\Di^2)}&\ds\le C\ \| \lf(1+|\nabla u|^2\rg)^{p/2-1}\ A\,\nabla u\|_{L^{\frac{p}{p-1}}(\Di^2)}\\[5mm]
\ds &\ds\le C\ \|  |\nabla u|+|\nabla u|^{p-1}\|_{L^{\frac{p}{p-1}}(\Di^2)}\\[5mm]
 &\ds\le C\ \lf[\|\nabla u\|_{L^{\frac{p}{p-1}}(\Di^2)}+\|\nabla u\|_{L^{p}(\Di^2)}^{p-1}  \rg]\ ,
\end{array}
\ee
where $C>0$ is independent of $p\in [2,3]$. Inserting this inequality in (\ref{IV.22}) gives also
\be
\label{IV.23-b}
\begin{array}{rl}
\ds\lf\|\lf(1+\|\nabla u\|_{L^p(\Di^2)}^2\rg)^{p/2-1}\,\nabla v\rg\|_{L^{\frac{p}{p-1}}(\Di^2)}&\le C\ \lf[\|\nabla u\|_{L^{\frac{p}{p-1}}(\Di^2)}+\|\nabla u\|_{L^{p}(\Di^2)}^{p-1}  \rg]\ .
\end{array}
\ee
From (\ref{IV.21}) we have
\be
\label{IV.25}
\||\nabla w|^{\frac{1}{p-1}}\|_{L^p(\Di^2)}\le\ C\  (p-2)^{\frac{1}{p-1}}\ \lf[ \|\nabla u\|_{L^p(\Di^2)}^{3/4}+ \|\nabla u\|_{L^p(\Di^2)}\rg]^{\frac{1}{p-1}}\ .
\ee
Using lemma~\ref{lm-holder-lorentz} for $g:=|\nabla w|^{\frac{1}{p-1}}$ gives
\be
\label{IV.26}
\||\nabla w|^{\frac{1}{p-1}}\|_{L^{2,1}(\Di^2)}\le\ C \ \pi^{\frac{p-2}{2p}}\ \lf( \frac{2\,(p-1)}{p-2}\rg)^{\frac{p-1}{p}}\ (p-2)^{\frac{1}{p-1}}\ \lf[ \|\nabla u\|_{L^p(\Di^2)}^{3/4}+ \|\nabla u\|_{L^p(\Di^2)}\rg]^{\frac{1}{p-1}}\ .
\ee
We have
\be
\label{IV.27}
\begin{array}{l}
\ds \frac{(p-2)^{\frac{1}{p-1}}}{(p-2)^{\frac{p-1}{p}}}=(p-2)^{\frac{1}{p-1}+\frac{1}{p}-1}\ .
\end{array}
\ee
Assume now $2<p<7/3$. This gives
\[
\frac{1}{p-1}+\frac{1}{p}>\frac{3}{4}+\frac{3}{7}=\frac{33}{28}>1\ .
\]
Hence there exists $\al>0$ independent of $2<p<7/3$ such that
\be
\label{IV.28}
\||\nabla w|^{\frac{1}{p-1}}\|_{L^{2,1}(\Di^2)}\le\ C \ (p-2)^\al\ \lf[ \|\nabla u\|_{L^p(\Di^2)}^{3/4}+ \|\nabla u\|_{L^p(\Di^2)}\rg]^{\frac{1}{p-1}}\ .
\ee
Observe that
\be
\label{IV.29}
\begin{array}{l}
|\nabla u|\le\lf| \lf(1+|\nabla u|^2\rg)^{p/2-1}\,|\nabla u| \rg|^{\frac{1}{p-1}}\ .
\end{array}
\ee
Hence from \eqref{IV.22} it follows that
\be\label{IV.30bis}
|\nabla u|\le \lf[ \lf(1+\|\nabla u\|_{L^p(\Di^2)}^2\rg)^{p/2-1}\rg]^{\frac{1}{p-1}}\lf(|\nabla v|+|A^{-1}\nabla^{\perp} w|\rg)^{\frac{1}{p-1}}\ee and
 for any $\om\subset\subset \Di^2$ we have
\be
\label{IV.30}
\|\nabla u\|_{L^{2,1}(\om)}\le \lf[ \lf(1+\|\nabla u\|_{L^p(\Di^2)}^2\rg)^{p/2-1}  \rg]^{\frac{1}{p-1}}\ \lf[\| |\nabla v|^\frac{1}{p-1}\|_{L^{2,1}(\om)}+  \| |\nabla w|^\frac{1}{p-1}\|_{L^{2,1}(\om)} \rg]\ .
\ee
By combining (\ref{IV.22}) and (\ref{IV.10})  gives
\be
\label{IV.31}
\begin{array}{l}
\ds  A\,\nabla\varphi+\nabla^\perp\psi= A\,\nabla\lf(\lf(1+\|\nabla u\|_{L^p(\Di^2)}^2\rg)^{p/2-1}\, v\rg)+ \nabla^\perp\lf(\lf(1+\|\nabla u\|_{L^p(\Di^2)}^2\rg)^{p/2-1}\, w\rg)\ .
\end{array}
\ee
We clearly have
\be
\label{IV.32}
\varphi-\lf(1+\|\nabla u\|_{L^p(\Di^2)}^2\rg)^{p/2-1}\, v \quad\mbox{ is $A$-harmonic }
\ee
and 
\be
\label{IV.32-a}
\psi-\lf(1+\|\nabla u\|_{L^p(\Di^2)}^2\rg)^{p/2-1}\, w\quad\mbox{ is $A^{-1}$-harmonic }\ .
\ee
This gives in particular using (\ref{IV.17-a}) and (\ref{IV.23-b})
\be
\label{IV.33}
\begin{array}{l}
\ds\lf(1+\|\nabla u\|_{L^p(\Di^2)}^2\rg)^{\frac{p}{2}-1}\ \|\nabla v\|_{L^{2,1}(\om)}\le \|\nabla \varphi\|_{L^{2,1}(\Di^2)}\\[5mm]
\ds\quad+C(\om)\, \lf\|\nabla\varphi-\lf(1+\|\nabla u\|_{L^p(\Di^2)}^2\rg)^{p/2-1}\,\nabla v\rg\|_{L^{\frac{p}{p-1}}(\Di^2)}\\[5mm]
\ds\quad\le  C(\om)\, \sqrt{\int_{\Di^2}\ (1+|\nabla u|^2)^{\frac{p}{2}-1}\ |\Om|^2\  dx^2}\ \ \sqrt{ \int_{\Di^2}\  \lf(1+|\nabla u|^2\rg)^{\frac{p}{2}}  dx^2} \\[5mm]
\ds\quad\ +\,C(\om)\ \lf[\|\nabla u\|_{L^{\frac{p}{p-1}}(\Di^2)}+\|\nabla u\|_{L^{p}(\Di^2)}^{p-1}  \rg]\ .
\end{array}
\ee
Observe that $p-1\ge 1$ and then $\beta:=1/(p-1)\le 1$. The decreasing rearangement of $|\nabla v|^\beta$ coincides obviously with the decreasing rearangement of $|\nabla v|$ to the power $\beta$ moreover one has obviously for any function $g$ on $[0,\pi]$
\be
\label{IV.33-a}
\int_0^{\pi}g(t)^\beta\frac{dt}{\sqrt{t}}\le \lf[ \int_0^{\pi}\frac{dt}{\sqrt{t}} \rg]^{1-\beta}\ \lf[ \int_0^{\pi}g(t)\frac{dt}{\sqrt{t}} \rg]^\beta\ .
\ee
Hence
\be
\label{IV.33-b}
\begin{array}{l}
 \ds\lf(1+\|\nabla u\|_{L^p(\Di^2)}^2\rg)^{\frac{p-2}{2(p-1)}}\ \||\nabla v|^\frac{1}{p-1}\|_{L^{2,1}(\om)}\\[5mm]
 \ds\quad\le C(\om)\ \lf[\int_{\Di^2}\ (1+|\nabla u|^2)^{\frac{p}{2}-1}\ |\Om|^2\  dx^2\rg]^{\frac{1}{2(p-1)}}\ \lf[\int_{\Di^2}\  \lf(1+|\nabla u|^2\rg)^{\frac{p}{2}}  dx^2\rg]^{\frac{1}{2(p-1)}}\\[5mm]
 \ds\quad+ C(\om)\ \lf[\|\nabla u\|^{\frac{1}{p-1}}_{L^{\frac{p}{p-1}}(\Di^2)}+\|\nabla u\|_{L^{p}(\Di^2)}  \rg]\ .
 \end{array}
\ee
This gives
\be
\label{IV.33-c}
\begin{array}{l}
 \ds\lf(1+\|\nabla u\|_{L^p(\Di^2)}^2\rg)^{\frac{p-2}{2(p-1)}}\ \||\nabla v|^\frac{1}{p-1}\|_{L^{2,1}(\om)}\\[5mm]
 \ds\quad\le C(\om)\ \lf[\int_{\Di^2}\ (1+|\nabla u|^2)^{\frac{p}{2}-1}\ |\Om|^2\  dx^2\rg]^{\frac{1}{2(p-1)}}\ \lf(1+\|\nabla u\|_{L^p(\Di^2)}^2\rg)^{\frac{p}{4(p-1)}}\\[5mm] 
 \ds \quad+\,C(\om)\ \|\nabla u\|_{L^{p}(\Di^2)}^{\frac{1}{p-1}}\ \lf(1+\|\nabla u\|_{L^p(\Di^2)}^2\rg)^{\frac{p-2}{2(p-1)}}\ ,
 \end{array}
 \ee
 \newpage
which gives
 \be
\label{IV.33-d}
\begin{array}{rl}
 \ds \||\nabla v|^\frac{1}{p-1}\|_{L^{2,1}(\om)}&\ds\le C(\om)\ \lf[\int_{\Di^2}\ (1+|\nabla u|^2)^{\frac{p}{2}-1}\ |\Om|^2\  dx^2\rg]^{\frac{1}{2(p-1)}}\ \lf(1+\|\nabla u\|_{L^p(\Di^2)}^2\rg)^{\frac{4-p}{4(p-1)}}\\[5mm] 
 \ds \quad\quad\quad&\ds+\,C(\om)\ \|\nabla u\|_{L^{p}(\Di^2)}^{\frac{1}{p-1}}\ .
 \end{array}
 \ee
 Combining (\ref{IV.28}), (\ref{IV.30}) and (\ref{IV.33-d}) gives
\be
\label{IV.34}
\begin{array}{l}
\ds\lf(1+\|\nabla u\|_{L^p(\Di^2)}^2\rg)^{\frac{2-p}{2(p-1)}}  \ \|\nabla u\|_{L^{2,1}(\om)}\le \ C \ (p-2)^\al\ \lf[ \|\nabla u\|_{L^p(\Di^2)}^{3/4}+ \|\nabla u\|_{L^p(\Di^2)}\rg]^{\frac{1}{p-1}}\\[5mm]
\ds+C(\om)\ \lf[\int_{\Di^2}\ (1+|\nabla u|^2)^{\frac{p}{2}-1}\ |\Om|^2\  dx^2\rg]^{\frac{1}{2(p-1)}}\ \lf(1+\|\nabla u\|_{L^p(\Di^2)}^2\rg)^{\frac{4-p}{4(p-1)}}+\,C(\om)\ \|\nabla u\|_{L^{p}(\Di^2)}^{\frac{1}{p-1}}\ .
\end{array}
\ee
This finally implies
\be
\label{IV.35}
\begin{array}{l}
\|\nabla u\|_{L^{2,1}(\om)}\le \  C \ (p-2)^\al\ \lf[ \|\nabla u\|_{L^p(\Di^2)}^{3/4}+ \|\nabla u\|_{L^p(\Di^2)}\rg]^{\frac{1}{p-1}}\  \lf(1+\|\nabla u\|_{L^p(\Di^2)}^2\rg)^{\frac{p-2}{2(p-1)}} \\[5mm]
\ds\ +\ C(\om)\ \lf[\int_{\Di^2}\ (1+|\nabla u|^2)^{\frac{p}{2}-1}\ |\Om|^2\  dx^2\rg]^{\frac{1}{2(p-1)}}\ \lf(1+\|\nabla u\|_{L^p(\Di^2)}^2\rg)^{\frac{p}{4(p-1)}} \\[5mm]
\ds\ +\,C(\om)\ \|\nabla u\|_{L^{p}(\Di^2)}^{\frac{1}{p-1}}\ \lf(1+\|\nabla u\|_{L^p(\Di^2)}^2\rg)^{\frac{p-2}{2(p-1)}}\ . 
\end{array}
\ee
We can conclude the {\bf proof of Lemma \ref{th-IV.1}.}~~~$\Box$
\subsection{The Case of Degenerating Annuli}
\begin{Lm}
\label{lm-IV.1-ann}
There exists $\sigma>0$ such that for any $p> 2$ the following holds. Let $0<\delta<1$ and denote
\be
\label{an-1}
(p-2)\,\log\delta^{-1}:=K
\ee
Let $\Om\in L^2(\Di^2,{\R}^2\otimes so(m))$ and $u\in W^{1,p}(\Di^2,{\R}^m)$ satisfying
\be
\label{I-r.1-ann}
-\,\mbox{div}\lf((1+|\nabla u|^2)^{\frac{p}{2}-1}\,\nabla u\rg)=(1+|\nabla u|^2)^{\frac{p}{2}-1}\,\Om\cdot \nabla u\ \mbox{ in }{\mathcal D}'(\Di^2)\ .
\ee
Assume
 \be
\label{I-r.2-ann}
\int_{B_1(0)\setminus B_\delta(0)}(1+|\nabla u|^2)^{\frac{p}{2}-1}\ |\Om|^2\ dx^2<\sigma\ ,
\ee
then for any $0<t<1$
\be
\label{an-2}
\lf\|  \lf[\sqrt{(1+|\nabla u|^2)^{\frac{p}{2}-1}}\,\lf|\frac{1}{r}\frac{\p u}{\p\theta}\rg|\rg]^\frac{1}{p-1}\rg\|_{L^{2,1}(B_t(0)\setminus B_{t^{-1}\,\delta}(0))}\le C(t,\|\nabla u\|_p,K)\ .
\ee
\hfill $\Box$
\end{Lm}
\noindent{\bf Proof of lemma~\ref{lm-IV.1-ann}.}
For any $1>\eta>0$ and $\eta^2>2\,\delta>0$ we denote
\be
\label{IV.36}
A(\eta,\delta)=B_\eta(0)\setminus B_{\delta/\eta}(0)\ .
\ee
We clearly concentrate to the case where $p$ is close to $2$ otherwise the lemma is trivial. Let then $p_k>2$ with $\lim_{k\rightarrow +\infty}p_k=2$ and let $\delta_k\rightarrow 0$. We shall omit to mention  the subscript $k$ unless it is necessary to evacuate any ambiguity. Let $\Om\in L^2(A(1,\delta),{\R}^2\otimes so(m))$ and $u\in W^{1,p}(D_2,{\R}^m)$ satisfying
\be
\label{IV.37}
-\,\mbox{div}\lf((1+|\nabla u|^2)^{\frac{p}{2}-1}\,\nabla u\rg)=(1+|\nabla u|^2)^{\frac{p}{2}-1}\,\Om\cdot \nabla u\ \quad\mbox{ in }{\mathcal D}'(A(1,\delta))\ .
\ee
Assume
 \be
\label{IV.38}
\int_{A(1,\delta_k)}(1+|\nabla u|^2)^{\frac{p}{2}-1}\ |\Om|^2\ dx^2<\sigma\ ,
\ee
where $\sigma>0$ is given by theorem~\ref{th-I.1}. We extend $\Om$ by zero inside $B_{\delta_k}(0)$ and we keep denoting $\Om$  this extension of $\Om$ by 0. Thanks to theorem~\ref{th-I.1} we have the existence of $(A,B)$ such that
\be
\label{IV.39}
(1+|\nabla u|^2)^{\frac{p}{2}-1}\,\lf(\nabla A+\,A\,\Om\rg)=\nabla^\perp B\quad\mbox{ in }{\mathcal D}'(\Di^2)\ ,
\ee
and
\be
\label{IV.40}
\lf\|\mbox{dist}(A,SO(m))\rg\|^2_{L^\infty(\Di^2)}+\int_{\Di^2}(1+|\nabla u|^2)^{\frac{p}{2}-1}\,|\nabla A|^2\ dx^2\le C\, \int_{\Di^2}(1+|\nabla u|^2)^{\frac{p}{2}-1}\ |\Om|^2\ dx^2\ ,
\ee
such that
\be
\label{IV.41}
\int_{\Di^2}\frac{|\nabla B|^2}{(1+|\nabla u|^2)^{\frac{p}{2}-1}\,}\ dx^2\le C\, \int_{\Di^2}(1+|\nabla u|^2)^{\frac{p}{2}-1}\ |\Om|^2\ dx^2\ .
\ee
Hence, because of (\ref{IV.39}) there holds
\be
\label{IV.42}
-\,\mbox{div}\lf((1+|\nabla u|^2)^{\frac{p}{2}-1}\,A\,\nabla u\rg)=\nabla^\perp B\cdot \nabla u\ \quad\mbox{ in }{\mathcal D}'(A(1,\delta))\ .
\ee
Let $\varphi$ be the solution of
\be
\label{IV.43}
\lf\{
\begin{array}{rl}
\ds -\mbox{div}(A\,\nabla\varphi)&\ds=\nabla^\perp B\cdot\nabla u\quad\mbox{ in }\Di^2\\[5mm]
\ds\varphi&\ds=0\quad\mbox{ on }\p \Di^2\ .
\end{array}
\rg.
\ee
As for (\ref{IV.12}) and (\ref{IV.12-a}) we have
\be
\label{IV.44}
\|\nabla B\|_{L^{p'}(\Di^2)}\le C_p \ \sqrt{\int_{\Di^2}\ (1+|\nabla u|^2)^{\frac{p}{2}-1}\ |\Om|^2\  dx^2}\ \ \lf[  \int_{\Di^2}\  \lf(1+|\nabla u|^2\rg)^{\frac{p}{2}}  dx^2 \rg]^{\frac{p-2}{2p}}\ .
\ee
and as (\ref{IV.17-a}) we deduce
\be
\label{IV.45}
\|\nabla\varphi\|_{L^{2,1}(\Di^2)}\le C\, \sqrt{\int_{\Di^2}\ (1+|\nabla u|^2)^{\frac{p}{2}-1}\ |\Om|^2\  dx^2}\ \ \sqrt{ \int_{\Di^2}\  \lf(1+|\nabla u|^2\rg)^{\frac{p}{2}}  dx^2} \ .
\ee
Because of (\ref{IV.42}) there exists $\psi$ on $A(1,\delta)$ such that
\be
\label{IV.46}
\nabla^\perp\psi=(1+|\nabla u|^2)^{\frac{p}{2}-1}\,A\,\nabla u-A\,\nabla\varphi- C_\ast\, \nabla\log r\quad\mbox{ in }{\mathcal D}'(A(1,\delta))\ .
\ee
where
\be
\label{IV.46-a}
2\pi\,C_\ast=\int_{\p B_r}(1+|\nabla u|^2)^{\frac{p}{2}-1}\,A\,\frac{\p u}{\p r}-A\,\frac{\p \varphi}{\p r} \ dl\ .
\ee
Observe that
\be
\label{IV.46-b}
\begin{array}{l}
\ds\int_{t^{-1}\,s^{-1}\,\delta}^{t\,s} \lf|\int_{\p B_r(0)}\,A\frac{\p\varphi}{\p r}\ dl\rg|^2\ r^{-1}\,dr\le \int_{t^{-1}\,s^{-1}\,\delta}^{t\,s} \int_{\p B_r(0)}\,\lf|A\frac{\p\varphi}{\p r}\rg|^2\ dl\,dr \\[5mm]
\ds\quad\le \int_{B_1(0)}|\nabla\varphi|^2\ dx^2\ .
\end{array}
\ee
and that
\be
\label{IV.46-c}
\begin{array}{l}
\ds\int_{t^{-1}\,s^{-1}\,\delta}^{t\,s}  \lf|\int_{\p B_r(0)}\,A\frac{\p\varphi}{\p r}\ dl\rg|\ r^{-1}\ dr\le \int_{ B_{ts}(0)\setminus B_{\delta/ts}(0)}\,\lf|\nabla\varphi\rg|(x) \ |x|^{-1}\ dx^2\\[5mm]
\ds\quad\le C\ \|\nabla \varphi\|_{L^{2,1}(\Di^2)}\ .
\end{array}
\ee
Hence we bound
\be
\label{IV.46-d}
\begin{array}{l}
\ds(2\pi)^{\frac{p}{p-1}}\,\int_{t^{-1}\,s^{-1}\,\delta}^{t\,s}\frac{|C_\ast|^{\frac{p}{p-1}}}{r}\ dr\\[5mm]
\ds\quad\le \int_{t^{-1}\,s^{-1}\,\delta}^{t\,s} \frac{1}{r}\ \lf[ \int_{\p B_r(0)}\,(1+|\nabla u|^2)^{\frac{p}{2}-1}\,\lf|A\,\nabla u\rg|+\lf|A\frac{\p\varphi}{\p r}\rg|\ dl\rg]^{\frac{p}{p-1}}\ dr \ .
\end{array}
\ee
Since $p/(p-1)<2$ we have for any $a>0,$ $a^{\frac{p}{p-1}}\le a+a^2$. This implies using (\ref{IV.46-b}) and (\ref{IV.46-c})
\be
\label{IV.46-e}
\begin{array}{l}
\ds\int_{t^{-1}\,s^{-1}\,\delta}^{t\,s} \frac{1}{r}\ \lf[ \int_{\p B_r(0)}\,\lf|A\frac{\p\varphi}{\p r}\rg|\ dl\rg]^{\frac{p}{p-1}}\ dr \\[5mm]
\ds\quad\le \int_{t^{-1}\,s^{-1}\,\delta}^{t\,s} \frac{1}{r}\  \int_{\p B_r(0)}\,\lf|A\frac{\p\varphi}{\p r}\rg|\ dl\ dr +\int_{t^{-1}\,s^{-1}\,\delta}^{t\,s} \frac{1}{r}\ \lf[ \int_{\p B_r(0)}\,\lf|A\frac{\p\varphi}{\p r}\rg|\ dl\rg]^{2}\ dr \\[5mm]
\ds\le C\,\lf [ \|\nabla\varphi\|_{L^{2,1}(\Di^2)}+\|\nabla\varphi\|^2_{L^{2}(\Di^2)}\rg]\ .
\end{array}
\ee
We have also
\be
\label{IV.46-f}
\begin{array}{l}
\ds\int_{t^{-1}\,s^{-1}\,\delta}^{t\,s} \frac{1}{r}\ \lf[ \int_{\p B_r(0)}\,(1+|\nabla u|^2)^{\frac{p}{2}-1}\,\lf|A\,\nabla u\rg|\ dl\rg]^{\frac{p}{p-1}}\ dr \\[5mm]
\ds\quad\le C\,\int_{t^{-1}\,s^{-1}\,\delta}^{t\,s} \frac{1}{r}\ \lf[ r^{\frac{1}{p}}\ \lf[\int_{\p B_r(0)}(1+|\nabla u|^2)^{\frac{p}{2}}\ dl\rg]^{\frac{p-1}{p}}\rg]^{\frac{p}{p-1}}\ dr\\[5mm]
\ds\quad\le C\,\int_{t^{-1}\,s^{-1}\,\delta}^{t\,s} r^{\frac{2-p}{p-1}}\ \int_{\p B_r(0)}(1+|\nabla u|^2)^{\frac{p}{2}}\ dl\ dr\\[5mm]
\ds\quad\le C\,\delta^{\frac{2-p}{p-1}}\,\int_{\Di^2}(1+|\nabla u|^2)^{\frac{p}{2}}\ dx^2=O(1)\ ,
\end{array}
\ee
where we used $\delta^{{2-p}}=O(1)$. This gives
\be
\label{IV.46-g}
|C_\ast|=O\lf(\frac{1}{\log^{\frac{p-1}{p}} \delta^{-1}}\rg)\ . 
\ee
In particular
\be
\label{IV.46-h}
\int_{B_1(0)\setminus B_\delta}\frac{|C_\ast|^2}{|x|^2}\ dx^2\le \log\delta^{-1}\ O\lf(\frac{1}{\log^{\frac{2p-2}{p}} \delta^{-1}}\rg) = O\lf(\frac{1}{\log^{\frac{p-2}{p}} \delta^{-1}}\rg) \ .
\ee
We have
\be
\label{IV.46-i}
\log^{\frac{p-2}{p}} \delta^{-1}=e^{\frac{p-2}{p}\, \log\log\delta^{-1}}=e^{O\lf(\frac{\log\log\delta^{-1}}{\log\delta^{-1}}\rg)}\rightarrow 1.
\ee
this finally gives that
\be
\label{lV.46-j}
\|\nabla(C_\ast\log|x|)\|_{L^2(B_1(0)\setminus B_\delta(0))}=O(1)\ .
\ee

Similarly as (\ref{IV.18})...(\ref{IV.30}) there exists  $v$ and $w$ in $W^{1,p'}(\Di^2,{\R}^m)$ such that
\be
\label{IV.51}
\frac{\ds \lf(1+|\nabla u|^2\rg)^{p/2-1}}{\ds\lf(1+\|\nabla u\|_{L^p(\Di^2)}^2\rg)^{p/2-1}}\ A\,\nabla u=A\,\nabla v+\nabla^\perp w\ ,
\ee
and
\be
\label{IV.52}
\||\nabla w|^{\frac{1}{p-1}}\|_{L^{2,1}(\Di^2)}\le\ C \ (p-2)^\al\ \lf[ \|\nabla u\|_{L^p(\Di^2)}^{3/4}+ \|\nabla u\|_{L^p(\Di^2)}\rg]^{\frac{1}{p-1}}
\ee
as well as\ .
Classical elliptic estimates applied to (\ref{IV.20}) give
\be
\label{IV.52-a}
\begin{array}{rl}
\ds\lf\|\lf(1+\|\nabla u\|_{L^p(\Di^2)}^2\rg)^{p/2-1}\,\nabla w\rg\|_{L^{\frac{p}{p-1}}(\Di^2)}&\ds\le C\ \| \lf(1+|\nabla u|^2\rg)^{p/2-1}\ A\,\nabla u\|_{L^{\frac{p}{p-1}}(\Di^2)}\\[5mm]
\ds &\ds\le C\ \|  |\nabla u|+|\nabla u|^{p-1}\|_{L^{\frac{p}{p-1}}(\Di^2)}\\[5mm]
 &\ds\le C\ \lf[\|\nabla u\|_{L^{\frac{p}{p-1}}(\Di^2)}+\|\nabla u\|_{L^{p}(\Di^2)}^{p-1}  \rg]\ ,
\end{array}
\ee
where $C>0$ is independent of $p\in [2,3]$. Inserting this inequality in (\ref{IV.51}) gives also
\be
\label{IV.52-b}
\begin{array}{rl}
\ds\lf\|\lf(1+\|\nabla u\|_{L^p(\Di^2)}^2\rg)^{p/2-1}\,\nabla v\rg\|_{L^{\frac{p}{p-1}}(\Di^2)}&\le C\ \lf[\|\nabla u\|_{L^{\frac{p}{p-1}}(\Di^2)}+\|\nabla u\|_{L^{p}(\Di^2)}^{p-1}  \rg]\ .
\end{array}
\ee

We have
\be
\label{IV.53}
\begin{array}{rl}
\ds\mbox{div}(A\,\nabla v)&\ds= A \,\mbox{div}\lf(\frac{\ds \lf(1+|\nabla u|^2\rg)^{p/2-1}}{\ds\lf(1+\|\nabla u\|_{L^p(\Di^2)}^2\rg)^{p/2-1}}\ \nabla u\rg)\\[5mm]
\ds\quad\quad\quad &\ds\quad+\frac{\ds \lf(1+|\nabla u|^2\rg)^{p/2-1}}{\ds\lf(1+\|\nabla u\|_{L^p(\Di^2)}^2\rg)^{p/2-1}}\ \nabla A\cdot\nabla u\quad\quad
\mbox{ on }{\mathcal D}'(\Di^2)\ .
\end{array}
\ee
Hence
\be
\label{IV.54}
\begin{array}{rl}
\ds\mbox{div}(A\,\nabla v)&\ds= \,\frac{\ds A\,\Om\cdot\nabla u}{\ds\lf(1+\|\nabla u\|_{L^p(\Di^2)}^2\rg)^{p/2-1}}\\[5mm]
\ds\quad\quad\quad &\ds\quad+\frac{\ds \lf(1+|\nabla u|^2\rg)^{p/2-1}}{\ds\lf(1+\|\nabla u\|_{L^p(\Di^2)}^2\rg)^{p/2-1}}\ \nabla A\cdot\nabla u\quad\quad
\mbox{ on }{\mathcal D}'(\Di^2)\ .
\end{array}
\ee
We have
\be
\label{IV.55}
\begin{array}{l}
\ds \int_{\Di^2} \lf| \lf(1+|\nabla u|^2\rg)^{p/2-1} \ \nabla A\cdot\nabla u\rg|\ dx^2\\[5mm]
\ds\quad\le \lf[ \int_{\Di^2} \lf(1+|\nabla u|^2\rg)^{p/2-1} \ |\nabla A|^2 \ dx^2\rg]^{\frac{1}{2}}\ \lf[ \int_{\Di^2} \lf(1+|\nabla u|^2\rg)^{p/2-1} \ |\nabla u|^2 \ dx^2\rg]^{\frac{1}{2}}\\[5mm]
\ds\quad\le \lf[\int_{B_1(0)\setminus B_\delta(0)} \lf(1+|\nabla u|^2\rg)^{p/2-1} \ |\Om|^2 \ dx^2\rg]^{\frac{1}{2}}\ \lf[ \int_{\Di^2} \lf(1+|\nabla u|^2\rg)^{p/2-1} \ |\nabla u|^2 \ dx^2\rg]^{\frac{1}{2}}\ .
\ds
\end{array}
\ee
This gives finally
\be
\label{IV.56}
\begin{array}{l}
\ds\lf(1+\|\nabla u\|_{L^p(\Di^2)}^2\rg)^{p/2-1}\ \|\mbox{div}(A\,\nabla v)\|_{L^1(\Di^2)}\\[5mm]
\ds\quad\le  \lf[\int_{\Di^2} \lf(1+|\nabla u|^2\rg)^{p/2-1} \ |\Om|^2 \ dx^2\rg]^{\frac{1}{2}}\ \lf[ \int_{\Di^2} \lf(1+|\nabla u|^2\rg)^{p/2-1} \ |\nabla u|^2 \ dx^2\rg]^{\frac{1}{2}}\\[5mm]
\end{array}
\ee
Let $\zeta$ be the solution of
\be
\label{IV.57}
\lf\{
\begin{array}{l}
\ds\Delta\zeta=\mbox{div}(A\,\nabla v)\quad\mbox{ in }{\mathcal D}'(\Di^2)\\[5mm] 
\ds\zeta=0\quad\mbox{ on }\p \Di^2\ .
\end{array}
\rg.
\ee
Classical elliptic estimates give 
\be
\label{IV.58}
\begin{array}{l}
\ds\lf(1+\|\nabla u\|_{L^p(\Di^2)}^2\rg)^{p/2-1}\ \|\nabla \zeta\|_{L^{2,\infty}(\Di^2)}\le C\ \lf(1+\|\nabla u\|_{L^p(\Di^2)}^2\rg)^{p/2-1}\ \|\Delta \zeta\|_{L^{1}(\Di^2)}\\[5mm]
\ds\quad\le C\ \lf[\int_{\Di^2} \lf(1+|\nabla u|^2\rg)^{p/2-1} \ |\Om|^2 \ dx^2\rg]^{\frac{1}{2}}\ \lf[ \int_{\Di^2} \lf(1+|\nabla u|^2\rg)^{p/2-1} \ |\nabla u|^2 \ dx^2\rg]^{\frac{1}{2}}\ .
\end{array}
\ee
Let $\xi$ be such that
\be
\label{IV.59}
A\,\nabla v= \nabla \zeta+\nabla^\perp\xi\ \quad\mbox{ in } \Di^2\ .
\ee
Because of (\ref{IV.52-b}) and (\ref{IV.58}) we have
\be
\label{IV.60}
\begin{array}{l}
\ds\lf(1+\|\nabla u\|_{L^p(\Di^2)}^2\rg)^{p/2-1}\ \|\nabla\xi\|_{L^{\frac{p}{p-1}}(\Di^2)}\\[5mm]
\ds\quad\le C\,  \lf[\int_{\Di^2} \lf(1+|\nabla u|^2\rg)^{p/2-1} \ |\Om|^2 \ dx^2\rg]^{\frac{1}{2}}\ \lf[ \int_{\Di^2} \lf(1+|\nabla u|^2\rg)^{p/2-1} \ |\nabla u|^2 \ dx^2\rg]^{\frac{1}{2}}\\[5mm]
\ds\quad+\,C\ \lf[\|\nabla u\|_{L^{\frac{p}{p-1}}(\Di^2)}+\|\nabla u\|_{L^{p}(\Di^2)}^{p-1}  \rg]\ .
\end{array}
\ee
The following equation holds
\be
\label{IV.61}
-\,\mbox{div}(A^{-1}\,\nabla \xi)=\nabla A^{-1}\cdot\nabla^\perp\zeta\quad\mbox{ in }{\mathcal D}'(\Di^2)\ .
\ee
Let $\eta$ be the solution given by lemma~\ref{lm-wente-type-estimates-bis} of
\be
\label{IV.62}
\lf\{
\begin{array}{l}
\ds -\,\mbox{div}(A^{-1}\,\nabla \eta)=\nabla A^{-1}\cdot\nabla^\perp\zeta\quad\mbox{ in }{\mathcal D}'(\Di^2)\\[5mm]
\ds\eta=0\quad\mbox{ on }\p \Di^2\ 
\end{array}
\rg.
\ee
and satisfy
\be
\label{IV.63}
\|\nabla \eta\|_{L^2(\Di^2)}\le C\ \|\nabla A^{-1}\|_{L^2(\Di^2)}\ \|\nabla\zeta\|_{L^{(2,\infty)}(\Di^2)}\ .
\ee
We have in particular since $C\ \|\nabla A^{-1}\|_{L^2(\Di^2)}$ is bounded by an universal constant
\be
\label{IV.64}
\begin{array}{l}
\ds\lf(1+\|\nabla u\|_{L^p(\Di^2)}^2\rg)^{p/2-1}\ \|\nabla(\xi-\eta)\|_{L^{\frac{p}{p-1}}(\Di^2)}\\[5mm]
\ds\quad\le C  \lf[\int_{\Di^2} \lf(1+|\nabla u|^2\rg)^{p/2-1} \ |\Om|^2 \ dx^2\rg]^{\frac{1}{2}}\ \lf[ \int_{\Di^2} \lf(1+|\nabla u|^2\rg)^{p/2-1} \ |\nabla u|^2 \ dx^2\rg]^{\frac{1}{2}}\\[5mm]
\ds\quad+\,C\ \lf[\|\nabla u\|_{L^{\frac{p}{p-1}}(\Di^2)}+\|\nabla u\|_{L^{p}(\Di^2)}^{p-1}  \rg]\ .
\end{array}
\ee
Combining this inequality with the homogeneous elliptic equation
\be
\label{IV.65}
-\,\mbox{div}(A^{-1}\,\nabla(\xi-\eta))=0\quad\mbox{ in }{\mathcal D}'(\Di^2)\ .
\ee
gives for any $t<1$ the existence of a constant  $C_t>0$ such that
\be
\label{IV.66}
\begin{array}{l}
\ds\lf(1+\|\nabla u\|_{L^p(\Di^2)}^2\rg)^{p/2-1}\ \|\nabla(\xi-\eta)\|_{L^{2}(B^2_t(0))}\\[5mm]
\ds\quad\le C_t  \lf[\int_{\Di^2} \lf(1+|\nabla u|^2\rg)^{p/2-1} \ |\Om|^2 \ dx^2\rg]^{\frac{1}{2}}\ \lf[ \int_{\Di^2} \lf(1+|\nabla u|^2\rg)^{p/2-1} \ |\nabla u|^2 \ dx^2\rg]^{\frac{1}{2}}\\[5mm]
\ds\quad+\,C_t\ \lf[\|\nabla u\|_{L^{\frac{p}{p-1}}(\Di^2)}+\|\nabla u\|_{L^{p}(\Di^2)}^{p-1}  \rg]\ .
\end{array}
\ee
Combining now (\ref{IV.58}), (\ref{IV.63}) and (\ref{IV.66}) is giving
\be
\label{IV.67}
\begin{array}{l}
\ds\lf(1+\|\nabla u\|_{L^p(\Di^2)}^2\rg)^{p/2-1}\ \|\nabla v\|_{L^{(2,\infty)}(B^2_t(0))}\\[5mm]
\ds\quad\le C_t  \lf[\int_{\Di^2} \lf(1+|\nabla u|^2\rg)^{p/2-1} \ |\Om|^2 \ dx^2\rg]^{\frac{1}{2}}\ \lf[ \int_{\Di^2} \lf(1+|\nabla u|^2\rg)^{p/2-1} \ |\nabla u|^2 \ dx^2\rg]^{\frac{1}{2}}\\[5mm]
\ds\quad+\,C_t\ \lf[\|\nabla u\|_{L^{\frac{p}{p-1}}(\Di^2)}+\|\nabla u\|_{L^{p}(\Di^2)}^{p-1}  \rg]\ .
\end{array}
\ee
Let
\be
\label{IV.67-a}
\hat{v}:=\lf(1+\|\nabla u\|_{L^p(\Di^2)}^2\rg)^{p/2-1}\ v\quad\mbox{ and }\quad\hat{w}:=\lf(1+\|\nabla u\|_{L^p(\Di^2)}^2\rg)^{p/2-1}\ w\ .
\ee
Hence we have
\be
\label{IV.68}
\mbox{div}(A\,\nabla(\hat{v}-\varphi))=0\quad\mbox{ in }A(1,\delta)
\ee
together with
\be
\label{IV.69}
\begin{array}{l}
\ds \|\nabla (\hat{v}-\varphi)\|_{L^{(2,\infty)}(B^2_t(0))}\\[5mm]
\ds\quad\le C_t  \lf[\int_{\Di^2} \lf(1+|\nabla u|^2\rg)^{p/2-1} \ |\Om|^2 \ dx^2\rg]^{\frac{1}{2}}\ \lf[ \int_{\Di^2} \lf(1+|\nabla u|^2\rg)^{p/2-1} \ |\nabla u|^2 \ dx^2\rg]^{\frac{1}{2}}\\[5mm]
\ds\quad+\,C_t\ \lf[\|\nabla u\|_{L^{\frac{p}{p-1}}(\Di^2)}+\|\nabla u\|_{L^{p}(\Di^2)}^{p-1}  \rg]\ .
\end{array}
\ee
Let $\la:=\psi- w $, it satisfies
\be
\label{IV.70}
A\, \nabla(\hat{v}-\varphi) - C_\ast\, \nabla\log r=\nabla^\perp \la\quad\mbox{ in }A(t,\delta)
\ee
We extend $\hat{v}$ outside $B^2_t(0)$ on the whole of ${\R}^2$ using a Whitney type extension $\ti{v}$ such that
\be
\label{IV.71}
\|\nabla \ti{v}\|_{L^{2,\infty}({\R}^2)}\le C\, \|\nabla \hat{v}\|_{L^{2,\infty}(B^2_t(0))}\ ,
\ee
and we extend $\varphi $ by 0 outside $\Di^2$. Finally we extend $A$ from $\Di^2$ to ${\R}^2$ using a Whitney extension in such a way that
\be
\label{IV.71-a}
\|\nabla\ti{A}\|_{L^2({\R}^2)}\le C\, \|\nabla{A}\|_{L^2(\Di^2)}\ .
\ee
This gives
\be
\label{IV.72}
\begin{array}{l}
\ds \|\nabla (\ti{v}-\ti{\varphi})\|_{L^{(2,\infty)}({\R}^2)}\\[5mm]
\ds\quad\le C_t  \lf[\int_{\Di^2} \lf(1+|\nabla u|^2\rg)^{p/2-1} \ |\Om|^2 \ dx^2\rg]^{\frac{1}{2}}\ \lf[ \int_{\Di^2} \lf(1+|\nabla u|^2\rg)^{p/2-1} \ |\nabla u|^2 \ dx^2\rg]^{\frac{1}{2}}\\[5mm]
\ds\quad+\,C_t\ \lf[\|\nabla u\|_{L^{\frac{p}{p-1}}(\Di^2)}+\|\nabla u\|_{L^{p}(\Di^2)}^{p-1}  \rg]\ .
\end{array}
\ee
Let $\mu$ be the solution such that $\nabla \mu\in L^2({\R}^2)$ of
\be
\label{IV.73}
-\,\Delta\mu = \nabla \ti{A}\cdot \nabla^\perp(\ti{v}-\ti{\varphi})\ .
\ee
We have
\be
\label{IV.74}
\|\nabla\mu\|_{L^{2}({\R}^2)}\le C\ \|\nabla A\|_{L^2(\Di^2)}\ \|\nabla (\hat{v}-{\varphi})\|_{L^{(2,\infty)}(B^2_t(0))}\ .
\ee
We have that $\la-\mu$ is harmonic on $A(t,\delta)$ and there holds
\be
\label{IV.75}
\|\nabla(\la-\mu)\|_{L^{2,\infty}(A(t,\delta))}\le C\ \|\nabla (\hat{v}-{\varphi})\|_{L^{(2,\infty)}(B^2_t(0))}+\|C_\ast\,\nabla\log|x|\|_{L^{2}(A(t,\delta))}\ .
\ee
Let $\la_0$ and $\mu_0$ be the zero Fourier modes of $\la$ and $\mu$ that is
\be
\label{IV.76}
\la_0(r):=\frac{1}{2\pi}\int_0^{2\pi} \la(r,\theta)\ d\theta  \quad\mbox{ and }\quad\mu_0(r):=\frac{1}{2\pi}\int_0^{2\pi} \mu(r,\theta)\ d\theta  \ .
\ee
Since $\la-\mu$ is harmonic there exists $C_0,C_1\in {\R}$ such that
\be
\label{IV.77}
\la_0-\mu_0(r)=C_0+C_1\,\log\,r\ .
\ee
Using (36) of \cite{LR} we have for any $s<1$
\be
\label{IV.78}
\begin{array}{l}
\ds\|\nabla (\la-\mu-\la_0-\mu_0)\|_{L^2(A(s\,t,\delta))}\le C_s \ \|\nabla (\la-\mu-\la_0-\mu_0)\|_{L^{2,\infty}(A(t,\delta))}\\[5mm]
\ds\quad\le C_s\ \|\nabla (\hat{v}-{\varphi})\|_{L^{(2,\infty)}(B^2_t(0))}+\|C_\ast\,\nabla\log|x|\|_{L^{2}(A(t,\delta))}\ .
\end{array}
\ee
Hence we deduce
\be
\label{IV.79}
\begin{array}{l}
\ds\|\nabla (\la-\la_0)\|_{L^2(A(s\,t,\delta))}\le C_s\ \|\nabla (\hat{v}-{\varphi})\|_{L^{(2,\infty)}(B^2_t(0))}+C\ \|\nabla A\|_{L^2(\Di^2)}\ 
\|\nabla (\hat{v}-{\varphi})\|_{L^{(2,\infty)}(B^2_t(0))}\\[5mm]
\ds\quad\quad+\|C_\ast\,\nabla\log|x|\|_{L^{2}(A(t,\delta))}\ .
\end{array}
\ee
We have for $r\in[t^{-1}\,s^{-1}\,\delta, t\,s]$
\be
\label{IV.80}
\frac{C_1}{r}=\dot{\la}_0(r)-\dot{\mu}_0(r)=-r^{-1}\,\int_0^{2\pi}\,A\, \p_\theta(\hat{v}-\varphi)\ d\theta-\dot{\mu}_0(r)=-r^{-1}\,\int_0^{2\pi}\,(A-\ov{A})\, \p_\theta(\hat{v}-\varphi)\ d\theta-\dot{\mu}_0(r)
\ee
where
\[
\ov{A}(r):=\dashint_0^{2\pi}A(r,\theta)\ d\theta\ .
\]
We have
\be
\label{IV.81}
\begin{array}{rl}
\ds r^{-1}\,\lf|\int_0^{2\pi}\,(A-\ov{A})\, \p_\theta(\hat{v}-\varphi)\ d\theta\rg|&\ds\le  \lf\|\frac{A-\ov{A}}{r}\rg\|_{L^\infty(\p B_r(0))}\ \int_{\p B_r} \lf| \nabla(\hat{v}-\varphi)\rg|\ dl\\[5mm]
\ds\ &\ds\le r^{-1}\ \int_{\p B_r} \lf| \nabla A\rg|\ dl\ \int_{\p B_r} \lf| \nabla(\hat{v}-\varphi)\rg|\ dl\ .
\end{array}
\ee
Thus
\be
\label{IV.82}
\forall\ t^{-1}\,s^{-1}\,\delta<r<t\,s\quad\quad |C_1|\  r^{-1}\le r^{-1}\ \int_{\p B_r} \lf| \nabla A\rg|\ dl\ \int_{\p B_r} \lf| \nabla(\hat{v}-\varphi)\rg|\ dl+|\dot{\mu}_0(r)|\ .
\ee
Choose $r\in [t^{-1}\,s^{-1}\,\delta, 2^{-1}\,t\,s]$, using the mean value theorem we obtain the existence of $\rho\in[r,2\,r]$ such that
\be
\label{IV.83}
\lf\{
\begin{array}{l}
\ds \int_{\p B_\rho} \lf| \nabla A\rg|\ dl\le C\,\lf[ \int_{B_{2r}(0)\setminus B_r(0)}|\nabla A|^2\ dx^2\rg]^{\frac{1}{2}}\\[5mm]
\ds\int_{\p B_r} \lf| \nabla(\hat{v}-\varphi)\rg|\ dl\le \|\nabla(\hat{v}-\varphi)\|_{L^{(2,\infty)}(B_2r(0)\setminus B_r(0))}\\[5mm]
\ds|\dot{\mu}_0(\rho)|\le C\ \dashint_r^{2r}|\dot{\mu}_0(r)|\ dr\ .
\end{array}
\rg.
\ee
where $C>0$ is a universal constant. Then we deduce
\be
\label{IV.84}
|C_1|\ \le C\,\lf[ \int_{B_{2r}(0)\setminus B_r(0)}|\nabla A|^2\ dx^2\rg]^{\frac{1}{2}}\  \|\nabla(\hat{v}-\varphi)\|_{L^{(2,\infty)}(B_t^2(0))}+\int_r^{2r}|\dot{\mu}_0(r)|\ dr\ .
\ee
This gives
\be
\label{IV.85}
|C_1|^2\ \le C\, \int_{B_{2r}(0)\setminus B_r(0)}|\nabla A|^2\ dx^2 \  \|\nabla(\hat{v}-\varphi)\|^2_{L^{(2,\infty)}(B_t^2(0))}+r\,\int_r^{2r}|\dot{\mu}_0(r)|^2\ dr\ .
\ee
Hence we deduce
\be
\label{IV.86}
|C_1|^2\ \log_2\lf(\frac{t^2\,s^2}{\delta}\rg) \le C\, \int_{A(ts,\delta)}|\nabla A|^2\ dx^2 \  \|\nabla(\hat{v}-\varphi)\|^2_{L^{(2,\infty)}(B_t^2(0))}+\int_{t^{-1}\,s^{-1}\,\delta}^{t\,s}|\dot{\mu}_0(r)|^2\ r\ dr\ .
\ee
This implies
\be
\label{IV.87}
\begin{array}{l}
\ds\int_{A(ts,\delta)}|\nabla(\la_0-\mu_0)|^2\ dx^2=\int_{A(ts,\delta)}|\nabla(C_1\,\log r)|^2\ dx^2\\[5mm]
\ds\quad\le C\, \|\nabla A\|^2_{L^2(B_1(0))}\  \|\nabla(\hat{v}-\varphi)\|^2_{L^{(2,\infty)}(B_t^2(0))}+\|\nabla \mu\|^2_{L^2(B_1(0))}\  .
\end{array}
\ee
Combining this last inequality with (\ref{IV.74}), (\ref{IV.75}) and (\ref{IV.78}) gives
\be
\label{IV.88}
\begin{array}{l}
\ds\|\nabla \la\|_{L^2(A(s\,t,\delta))}
\ds\le C_s\ \|\nabla (\hat{v}-{\varphi})\|_{L^{(2,\infty)}(B^2_t(0))}+C\ \|\nabla A\|_{L^2(\Di^2)}\ \|\nabla (\hat{v}-{\varphi})\|_{L^{(2,\infty)}(B^2_t(0))}\\[5mm]
\ds\quad+\|C_\ast\,\nabla\log|x|\|_{L^{2}(A(t,\delta))}\ .
\end{array}
\ee
Combining this fact together with (\ref{IV.45}) (\ref{IV.69}) and (\ref{lV.46-j}) is giving
\be
\label{IV.89}
\|\nabla \hat{v}\|_{L^2(A(s\,t,\delta))}=O(1)\ .
\ee
We construct a Whitney extension for $\hat{v}$ outside $A(s\,t,\delta)$ over the whole ${\R}^2$  that we denote $\check{v}$ such that
\be
\label{IV.90}
\|\nabla\check{v}\|_{L^{2}({\R}^2)}\le C\ \|\nabla \hat{v}\|_{L^2(A(s\,t,\delta))}\ .
\ee
Let $\nu$ be the solution such that $\nabla \nu\in L^2({\R}^2)$ of
\be
\label{IV.91}
-\,\Delta\nu = \nabla \ti{A}\cdot \nabla^\perp(\check{v}-\ti{\varphi})\ .
\ee
We have using \cite{CLMS}
\be
\label{IV.92}
\|\nabla\nu\|_{L^{2,1}({\R}^2)}\le\ \|\nabla A\|_{L^2(\Di^2)}\ \lf[\|\nabla\hat{v}\|_{L^2(A(s\,t,\delta))}+\|\nabla\varphi\|_{L^2(\Di^2)}\rg]\ .
\ee
We have that $\la-\nu$ is harmonic on $A(s\,t,\delta)$ and we denote similarly as before
\be
\label{IV.93}
\la_0(r):=\frac{1}{2\pi}\int_0^{2\pi} \la(r,\theta)\ d\theta  \quad\mbox{ and }\quad\nu_0(r):=\frac{1}{2\pi}\int_0^{2\pi} \nu(r,\theta)\ d\theta  \ .
\ee
Since $\la-\nu$ is harmonic there exists $C_0,C_1\in {\R}$ such that
\be
\label{IV.94}
\la_0-\nu_0(r)=\check{C}_0+\check{C}_1\,\log\,r\ .
\ee
Using (36) of \cite{LR} we have for any $\tau<s$
\be
\label{IV.95}
\begin{array}{l}
\ds\|\nabla (\la-\nu-\la_0-\nu_0)\|_{L^{2,1}(A(\tau\,t,\delta))}\le C_\tau \ \|\nabla (\la-\nu-\la_0-\nu_0)\|_{L^{2}(A(s\,t,\delta))}\\[5mm]
\ds\quad\le C_s\ \|\nabla (\hat{v}-{\varphi})\|_{L^{(2,\infty)}(B^2_t(0))}+\|C_\ast\,\nabla\log|x|\|_{L^{2}(A(t,\delta))}\ .
\end{array}
\ee
Hence we deduce 
\be
\label{IV.96}
\begin{array}{l}
\|\nabla (\la-\la_0)\|_{L^{2,1}(A(\tau\,t,\delta))}\le C_\tau\ \|\nabla (\hat{v}-{\varphi})\|_{L^{(2,\infty)}(B^2_t(0))}+C\ \|\nabla A\|_{L^2(\Di^2)}\ \|\nabla (\hat{v}-{\varphi})\|_{L^{(2,\infty)}(B^2_t(0))}\\[5mm]
\ds\quad\quad+\|C_\ast\,\nabla\log|x|\|_{L^{2}(A(t,\delta))}\ .
\end{array}
\ee

We have for $r\in[t^{-1}\,\tau^{-1}\,\delta, t\,\tau]$
\be
\label{IV.97}
\frac{\check{C}_1}{r}=\dot{\la}_0(r)-\dot{\nu}_0(r)=-r^{-1}\,\int_0^{2\pi}\,A\, \p_\theta(\hat{v}-\varphi)\ d\theta-\dot{\nu}_0(r)=-r^{-1}\,\int_0^{2\pi}\,(A-\ov{A})\, \p_\theta(\hat{v}-\varphi)\ d\theta-\dot{\nu}_0(r)\ ,
\ee
where
\[
\ov{A}(r):=\dashint_0^{2\pi}A(r,\theta)\ d\theta\ .
\]
We have
\be
\label{IV.98}
\begin{array}{rl}
\ds r^{-1}\,\lf|\int_0^{2\pi}\,(A-\ov{A})\, \p_\theta(\hat{v}-\varphi)\ d\theta\rg|&\ds\le  \lf\|\frac{A-\ov{A}}{r}\rg\|_{L^\infty(\p B_r(0))}\ \int_{\p B_r} \lf| \nabla(\hat{v}-\varphi)\rg|\ dl\\[5mm]
\ds\ &\ds\le r^{-1}\ \int_{\p B_r} \lf| \nabla A\rg|\ dl\ \int_{\p B_r} \lf| \nabla(\hat{v}-\varphi)\rg|\ dl\ .
\end{array}
\ee
Thus
\be
\label{IV.99}
\forall\ t^{-1}\,\tau^{-1}\,\delta<r<t\,\tau\quad\quad |\check{C}_1|\  \rho^{-1}\le \rho^{-1}\ \int_{\p B_\rho} \lf| \nabla A\rg|\ dl\ \int_{\p B_\rho} \lf| \nabla(\hat{v}-\varphi)\rg|\ dl+|\dot{\nu}_0(r)|\ .
\ee
Choose $r\in [t^{-1}\,\tau^{-1}\,\delta, 2^{-1}\,t\,\tau]$, using the mean value theorem we obtain the existence of $\rho\in[r,2\,r]$ such that
\be
\label{IV.100}
\lf\{
\begin{array}{l}
\ds \int_{\p B_\rho} \lf| \nabla A\rg|\ dl\le C\,\lf[ \int_{B_{2r}(0)\setminus B_r(0)}|\nabla A|^2\ dx^2\rg]^{\frac{1}{2}}\\[5mm]
\ds\int_{\p B_r} \lf| \nabla(\hat{v}-\varphi)\rg|\ dl\le \lf[ \int_{B_{2r}(0)\setminus B_r(0)}|\nabla(\hat{v}-\varphi)|^2\ dx^2\rg]^{\frac{1}{2}}\\[5mm]
\ds|\dot{\nu}_0(\rho)|\le C\ \dashint_r^{2r}|\dot{\nu}_0(r)|\ dr\ .
\end{array}
\rg.
\ee
where $C>0$ is a universal constant. Then we deduce
\be
\label{IV.101}
|\check{C}_1|\ \le C\,\lf[ \int_{B_{2r}(0)\setminus B_r(0)}|\nabla A|^2\ dx^2\rg]^{\frac{1}{2}}\  \lf[ \int_{B_{2r}(0)\setminus B_r(0)}|\nabla(\hat{v}-\varphi)|^2\ dx^2\rg]^{\frac{1}{2}}+\int_r^{2r}|\dot{\nu}_0(r)|\ dr\ .
\ee
This gives
\be
\label{IV.102}
\begin{array}{l}
\ds\int_{t^{-1}\,\tau^{-1}\,\delta}^{t\,\tau}\frac{|\check{C}_1|}{r}\ dr\le C\,\int_{t^{-1}\,\tau^{-1}\,\delta}^{t\,\tau}r^{-1}\, \int_{B_{2r}(0)\setminus B_r(0)}|\nabla A|^2\ dx^2\ dr\\[5mm]
\ds+C\,\int_{t^{-1}\,\tau^{-1}\,\delta}^{t\,\tau}r^{-1}\, \int_{B_{2r}(0)\setminus B_r(0)}|\nabla (\hat{v}-\varphi)|^2\ dx^2\ dr+\int_{t^{-1}\,\tau^{-1}\,\delta}^{t\,\tau}r^{-1}\int_r^{2r}|\dot{\nu}_0(s)|\ ds\ dr \ .
\end{array}
\ee
We deduce
\be
\label{IV.103}
|\check{C}_1|\,\log\lf( \frac{(t\,\tau)^2}{\delta} \rg)\le C\, \int_{\Di^2}|\nabla A|^2\ dx^2+\int_{A(s\,t,\delta)}|\nabla\hat{v}|^2\ dx^2+\int_{\Di^2}|\nabla\varphi|^2\ dx^2+\|\nabla \nu\|_{L^{2,1}({\R}^2)}\ .
\ee
This implies
\be
\label{IV.104}
\begin{array}{l}
\ds\|\nabla (\la_0-\nu_0)\|_{L^{2,1}(A(s\,t,\delta))}\\[5mm]
\ds\quad\le C\ \int_{\Di^2}|\nabla A|^2\ dx^2+\int_{A(s\,t,\delta)}|\nabla\hat{v}|^2\ dx^2+\int_{\Di^2}|\nabla\varphi|^2\ dx^2+\|\nabla \nu\|_{L^{2,1}({\R}^2)}\ .
\end{array}
\ee
Combining (\ref{IV.92}), (\ref{IV.96}) and (\ref{IV.104}) is giving
\be
\label{IV.105}
\begin{array}{l}
\|\nabla \la\|_{L^{2,1}(A(\tau\,t,\delta))}\le C_\tau\ \|\nabla (\hat{v}-{\varphi})\|_{L^{(2,\infty)}(B^2_t(0))}+C\ \|\nabla A\|_{L^2(\Di^2)}\ \|\nabla (\hat{v}-{\varphi})\|_{L^{(2,\infty)}(B^2_t(0))}\\[5mm]
\ds\quad+ C\ \int_{\Di^2}|\nabla A|^2\ dx^2+C\, \int_{A(s\,t,\delta)}|\nabla\hat{v}|^2\ dx^2+\int_{\Di^2}|\nabla\varphi|^2\ dx^2+\|C_\ast\,\nabla\log|x|\|_{L^{2}(A(t,\delta))}\ .
\end{array}
\ee
Combining this fact with the identity (\ref{IV.70}) is implying
\be
\label{IV.106}
\lf\|\frac{1}{r}\frac{\p \hat{v}}{\p \theta}\rg\|_{L^{2,1}(A(\tau\,t,\delta))}=O(1)\ .
\ee
Combining this fact with (\ref{IV.52}) is giving
\be
\label{IV.107}
\lf\| (1+|\nabla u|^2)^{\frac{p-2}{2p-2}}\lf|\frac{1}{r}\frac{\p u}{\p \theta}\rg|^{\frac{1}{p-1}}\rg\|_{L^{2,1}(A(\tau\,t,\delta))}=O(1)\ .
\ee
This concludes the proof of lemma~\ref{lm-IV.1-ann}.\hfill $\Box$
\section{Uniform $\epsilon $-regularity for $p $-harmonic $\Om $-systems as $p\rightarrow 2$.}
\reset
We give a consequence of the uniform $L^{2,1} $-bound on the disc (lemma~\ref{th-IV.1}).
\begin{Th}
\label{th-eps-reg}
Let $m\in {\N}$. Let $p\in[2,3]$ and $u\in W^{1,p}(\Di^2,{\R}^m)$ and denote
\be
\label{ep-0}
E:=\int_{\Di^2}(1+|\nabla u|^2)^{\frac{p}{2}-1}\ |\nabla u|^2\ dx^2\ .
\ee
There exists $\sigma>0$ depending only on $m$ and $E$ such that, for any  $\Om\in L^2(\Di^2,{\R}^2\otimes so(m))$ satisfying
\be
\label{ep-1}
-\,\mbox{div}\lf((1+|\nabla u|^2)^{\frac{p}{2}-1}\,\nabla u\rg)=(1+|\nabla u|^2)^{\frac{p}{2}-1}\,\Om\cdot \nabla u\ \mbox{ in }{\mathcal D}'(\Di^2)\ .
\ee
and the following conditions:     $|\Om|\le C\, |\nabla u|$ for some $C>0$, 
 \be
\label{ep-2}
\int_{\Di^2}(1+|\nabla u|^2)^{\frac{p}{2}-1}\ |\Om|^2\ dx^2\le \sigma\quad \mbox{ and }\quad\int_{\Di^2}|\nabla u|^2\ dx^2<\sigma\ .
\ee
Then
\be
\label{ep-02}
\|\nabla u\|^2_{L^\infty(B_{1/2}(0))}\le C(E,m)\ .
\ee
\hfill $\Box$
\end{Th}
\noindent{\bf Proof of theorem~\ref{th-eps-reg}.} Let $\varphi$ be the solution of
\be
\label{ep-3}
\lf\{
\begin{array}{l}
\ds\Delta\varphi=(1+|\nabla u|^2)^{\frac{p}{2}-1}\,\Om\cdot \nabla u\quad \mbox{ in }{\mathcal D}'(\Di^2)\\[5mm]
\ds\varphi=0\quad\mbox{ on }\p \Di^2\ .
\end{array}
\rg.
\ee
Classical elliptic estimates give
\be
\label{ep-4}
\begin{array}{l}
\ds\|\nabla \varphi\|_{L^{2,\infty}(\Di^2)}\le C\, \int_{\Di^2}|1+|\nabla u|^2|^{\frac{p}{2}-1}\, |\Om|\,|\nabla u|\ dx^2\\[5mm]
\ds\quad\le C\, \sqrt{\int_{\Di^2}(1+|\nabla u|^2)^{\frac{p}{2}-1}\ |\Om|^2\ dx^2}\  \sqrt{\int_{\Di^2}(1+|\nabla u|^2)^{\frac{p}{2}-1}\ |\nabla u|^2\ dx^2}\ .
\end{array}
\ee
By construction there holds
\be
\label{ep-5}
-\,\mbox{div}\lf((1+|\nabla u|^2)^{\frac{p}{2}-1}\,\nabla u\rg)=\mbox{div}\lf(\nabla \varphi\rg)\ .
\ee
Thanks to the main result in \cite{DM} we have for any $1\le q<2$
\be
\label{ep-6}
\begin{array}{l}
\ds\|\nabla u\|^{p-1}_{L^{q\,(p-1)}(B_{3/4}(0))}=\||\nabla u|^{p-1}\|_{L^q(B_{3/4}(0))}\le C\ \|\nabla \varphi\|_{L^q(\Di^2)}+\, C\,\|\nabla u\|_{L^2(\Di^2)}\\[5mm]
\ds\quad\le C_q\, \sqrt{\int_{\Di^2}(1+|\nabla u|^2)^{\frac{p}{2}-1}\ |\Om|^2\ dx^2}\  \sqrt{\int_{\Di^2}(1+|\nabla u|^2)^{\frac{p}{2}-1}\ |\nabla u|^2\ dx^2}+\, C\,\|\nabla u\|_{L^2(\Di^2)}\ .
\end{array}
\ee
We fix $q$ such that 
\[
p<p+\frac{p-2}{2}=q\,(p-1)<p+(p-2)=2\, (p-1)\ .
\]
Let $\delta>0$ such that
\be
\label{ep-7}
\frac{1}{p}=\frac{\delta}{q\,(p-1)}+\frac{1-\delta}{2}\ .
\ee
Explicit computations give
\be
\label{ep-8}
\delta=1-\frac{2}{3p}\ .
\ee
Classical H\"older inequality is giving
\be
\label{ep-9}
\begin{array}{l}
\ds\|\nabla u\|_{L^p(B_{3/4}(0))}\le\|\nabla u\|_{L^{q\,(p-1)}(B_{3/4}(0))}^\delta\ \|\nabla u\|_2^{1-\delta}\\[5mm]
\ds\le C\, \lf[ \int_{\Di^2}(1+|\nabla u|^2)^{\frac{p}{2}-1}\ |\Om|^2\ dx^2\rg]^{\frac{\delta}{2p-2}}\, \lf[\int_{\Di^2}(1+|\nabla u|^2)^{\frac{p}{2}-1}\ |\nabla u|^2\ dx^2\rg]^{\frac{\delta}{2p-2}}\ \|\nabla u\|_2^{1-\delta}\\[5mm]
\quad+C\, \|\nabla u\|_2^{1-\delta\,\frac{p-2}{p-1}}\le C\, (\sqrt{\sigma})^{1-\delta\,\frac{p-2}{p-1}}\ \lf[1+E^{\frac{\delta}{2p-2}}\rg]\ .
\end{array}
\ee
Combining this fact with (\ref{IV.35}) is giving
\be
\label{ep-10}
\|\nabla u\|_{L^{2,1}(B_{1/2}(0))}\le C(E)\ \sigma^\beta\ .
\ee
For some $\beta>0$ uniformly bounded from below for $p\in(2,3]$ and for $0<\sigma<1$. Let $x_0\in B_{1/4}(0)$ and $r<1/4$. Multiplying the equation by $\chi_{x_0,r}(x)\, (u-\ov{u}_{x_0,r}))$ where 
\[
\ov{u}_{x_0,r}:=\dashint_{B_r(x_0)\setminus B_{r/2}(x_0)} u(x)\ dx^2\ ,\ \chi_{x_0,r}:=\chi\lf(\frac{|x-x_0|}{r}\rg)
\]
and $\chi\in C^\infty_0([-1,1])$ with $\chi\ge 0$ and $\chi\equiv 1$ on $[-1/2,1/2]$ is giving
\be
\label{ep-11}
\begin{array}{l}
\ds\int_{B_r(x_0)}\chi_{x_0,r}(x)\ (1+|\nabla u|^2)^{\frac{p}{2}-1}\,|\nabla u|^2\, dx^2\\[5mm]
\ds\quad\le r^{-1}\, \|\chi'\|_\infty\ \int_{B_r(x_0)\setminus B_{r/2}(x_0)}\ |u-\ov{u}_{x_0,r}|\ (1+|\nabla u|^2)^{\frac{p}{2}-1}\,|\nabla u|\ dx^2\\[5mm]
\ds\quad+\int_{B_r(x_0)}\chi_{x_0,r}(x)\ \ |u-\ov{u}_{x_0,r}|\ (1+|\nabla u|^2)^{\frac{p}{2}-1}\ |\Om\cdot\nabla u|\ dx^2\ .
\end{array}
\ee
Hence using (\ref{ep-10}) and the injection of $W^{1,(2,1)}(\Di^2)$ into $C^0$ we obtain
\be
\label{ep-12}
\begin{array}{l}
\ds\int_{B_{r/2}(x_0)}(1+|\nabla u|^2)^{\frac{p}{2}-1}\,|\nabla u|^2\, dx^2\\[5mm]
\ds \le\ C\, \lf[\int_{B_r(0)\setminus B_{r/2}(0)}r^{-p}\,|u-\ov{u}_{x_0,r}|^p\ dx^2\rg]^{\frac{1}{p}}\  \lf[\int_{B_r(0)\setminus B_{r/2}(0)}\ \lf[|\nabla u|^{\frac{p}{p-1}}+|\nabla u|^p\rg] dx^2\rg]^{\frac{p-1}{p}}\\[5mm]
\ds+\ C(E)\ \sigma^\beta\ \int_{B_{r}(x_0)}(1+|\nabla u|^2)^{\frac{p}{2}-1}\,|\nabla u|^2\, dx^2\ .
\end{array}
\ee
Applying Poincar\'e inequality, for $\sigma^\beta$ sufficiently small we finally get
\be
\label{ep-13}
\int_{B_{r/2}(x_0)}(1+|\nabla u|^2)^{\frac{p}{2}-1}\,|\nabla u|^2\, dx^2\le C\ \int_{B_{r}(x_0)\setminus B_{r/2}(x_0) }(1+|\nabla u|^2)^{\frac{p}{2}-1}\,|\nabla u|^2\, dx^2\ .
\ee
This inequality implies a Morrey decrease of the form
\be
\label{ep-14}
\int_{B_r(x_0)} (1+|\nabla u|^2)^{\frac{p}{2}-1}\,|\nabla u|^2\, dx^2\le C\ r^\al
\ee
for some $\al$ independent of $p$ as $p\rightarrow 2$. This Morrey estimate is implying a uniform H\"older type estimate independent of $p\rightarrow 2$ and the result follows from standard bootstrap arguments.\hfill $\Box$
\section{Neck Length Control : Proof of Theorem~\ref{th-IV.2-ann}.}
\reset
We can assume that
\be
\label{IV.107-a}
\int_{B_1(0)\setminus B_\delta(0)}(1+|\nabla u|^2)^{\frac{p}{2}-1}\ |\Om|^2\ dx^2<\sigma\ .
\ee
otherwise we decompose the annulus $B_1(0)\setminus B_\delta(0)$ into $Q/\sigma$ annuli (as it is done in \cite{LR}).
We write
\be
\label{IV.108}
(1+|\nabla u|^2)^{\frac{p}{4}-\frac{1}{2}}\ \lf[(1+|\nabla u|^2)^{\frac{p}{2}-1}\,\lf|\frac{1}{r}\frac{\p u}{\p\theta}\rg|\rg]^\frac{1}{p-1}\ge\sqrt{(1+|\nabla u|^2)^{\frac{p}{2}-1}}\ \lf|\frac{1}{r}\frac{\p u}{\p\theta}\rg|\ .
\ee
Because of the hypothesis we have
\be
\label{IV.109}
(1+|\nabla u|^2)^{\frac{p}{4}-\frac{1}{2}}\le C\ \lf[1+\|\nabla u\|^{\frac{p-2}{2}}_\infty\rg]\le C\ \lf[1+L^{\frac{p-2}{2}}\ e^{\frac{p-2}{2} \,\log \delta^{-1}}\rg]\le C\ \lf[1+\sqrt{ L^{p-2}\ e^K}\rg]\ .
\ee
Combining (\ref{IV.107}) and (\ref{an-2}) is giving for any $0<t<1$
\be
\label{IV.1009}
\begin{array}{l}
\ds\lf\|(1+|\nabla u|^2)^{\frac{p}{4}-\frac{1}{2}}\,\lf|\frac{1}{r}\frac{\p u}{\p\theta}\rg|\rg\|_{L^{2,1}(B_t(0)\setminus B_{t^{-1}\,\delta}(0))}\le C(t,\|\nabla u\|_p,K, (p-2)\,\log L)\ .
\end{array}
\ee
Since we are assuming $u$ lipschitz and solving (\ref{II-r.1-ann}) and since $\Om$ is antisymmetric, the following Pohozaev type inequality holds as for the harmonic map case (see \cite{Lam})
\be
\label{IV.1100}
\begin{array}{l}
\ds\forall \,r\in[\delta,1]\quad\int_{\p B_r(0)}\, (1+|\nabla u|^2)^{\frac{p}{2}-1}\,\lf|\frac{\p u}{\p r}\rg|^2\ dl\le C\,\int_{\p B_r(0)}\, (1+|\nabla u|^2)^{\frac{p}{2}-1}\,\lf|\frac{1}{r}\frac{\p u}{\p\theta}\ \rg|^2\ dl\\[5mm]
\ds\quad\quad+C\ \frac{p-2}{r}\int_{B_r(0)}(1+|\nabla u|^2)^{\frac{p}{2}}\ dx^2\ .
\end{array}
\ee
Hence
\be
\label{IV.110}
\begin{array}{rl}
\ds\forall \,r\in[\delta,1]\quad\sqrt{\int_{\p B_r(0)}\, (1+|\nabla u|^2)^{\frac{p}{2}-1}\,\lf|\frac{\p u}{\p r}\rg|^2\ dl}&\ds\le C\,\sqrt{\int_{\p B_r(0)}\, (1+|\nabla u|^2)^{\frac{p}{2}-1}\,\lf|\frac{1}{r}\frac{\p u}{\p\theta}\ \rg|^2\ dl}\\[5mm]
\ds\quad\quad&\ds+C\ \frac{\sqrt{p-2}}{\sqrt{r}}\ .
\end{array}
\ee
Recall from (\ref{IV.46-a})
\be
\label{IV.46-rep}
2\pi\,C_\ast=\int_{\p B_r}(1+|\nabla u|^2)^{\frac{p}{2}-1}\,A\,\frac{\p u}{\p r}-A\,\frac{\p \varphi}{\p r} \ dl\ .
\ee
Hence in particular we deduce from (\ref{IV.110})
\be
\label{IV.111}
\begin{array}{l}
\ds\frac{|C_\ast|}{\sqrt{r}}\le C\,\sqrt{  \|(1+|\nabla u|^2)^{\frac{p}{2}-1}\|_\infty }\ \sqrt{\int_{\p B_r(0)}\, (1+|\nabla u|^2)^{\frac{p}{2}-1}\,\lf|\frac{\p u}{\p r}\rg|^2\ dl}+C\,\sqrt{\int_{\p B_r(0)}\, \lf|\frac{\p \varphi}{\p r}\rg|^2\ dl} \\[5mm]
\ds\le C\ \lf[1+\sqrt{ L^{p-2}\ e^K}\rg]\   \sqrt{\int_{\p B_r(0)}\, (1+|\nabla u|^2)^{\frac{p}{2}-1}\,\lf|\frac{1}{r}\frac{\p u}{\p\theta}\ \rg|^2\ dl}  \\[5mm]
\ds+C\ \lf[1+\sqrt{ L^{p-2}\ e^K}\rg] \frac{\sqrt{p-2}}{\sqrt{r}}+C\,\sqrt{\int_{\p B_r(0)}\, \lf|\frac{\p \varphi}{\p r}\rg|^2\ dl}\ .\end{array}
\ee
Using (\ref{IV.1009}), (\ref{IV.45}), lemma 9.2 from \cite{Mi} and the fact that
\be
\label{IV.0111}
\lf\| \sqrt{r^{-1}}  \rg\|_{L^{2,1}([\delta,1])}\simeq \log\frac{1}{\delta}\ ,
\ee
we deduce
\be
\label{IV.112}
|C_\ast|\, \log\delta^{-1}\, \le C( L,K, \|\nabla u\|_{L^p})\, \lf(1+\sqrt{p-2}\ \log\delta^{-1}\rg)\ .
\ee
Inserting this estimate in (\ref{IV.70}) and using (\ref{IV.105}) we obtain for any $0<t<1$
\be
\label{IV.113}
\|\nabla v\|_{L^{2,1}(A(t,\delta))}\le C( t,L,K, \|\nabla u\|_{L^p})\, \lf(1+\sqrt{p-2}\ \log\delta^{-1}\rg)\ .
\ee
Combining this fact with (\ref{IV.51}) and (\ref{IV.52} is finally giving
\be
\label{IV.114}
\begin{array}{l}
\ds\|\nabla u\|_{L^{2,1}(A(t,\delta))}\le \lf\| (1+|\nabla u|^2)^{\frac{p-2}{2p-2}}\lf|\nabla u\rg|^{\frac{1}{p-1}}\rg\|_{L^{2,1}(A(\tau\,t)}\\[5mm]
\ds\quad\le C( t,L,K, \|\nabla u\|_{L^p})\, \lf(1+\sqrt{p-2}\ \log\delta^{-1}\rg)\ .
\end{array}
\ee
Using proposition 2.7 of \cite{MiRi} we deduce (\ref{II.r.4-ann}) and theorem~\ref{th-IV.2-ann} is proved. \hfill $\Box$
\appendix
\addcontentsline{toc}{section}{Appendices}
\section*{Appendix}

\section{A Weighted Wente Inequality}
\reset
\begin{Lma}
\label{lm-wente-weight}
Let $p> 2$. Let $f$ be a non negative function on the disc $\Di^2$ satisfying $f\ge 1$. Assume $f\in L^{p/(p-1)}(\Di^2)$. Let $a$ and $b$ such that $\nabla a\in L^{p/(p-1)}(\Di^2)$ and $\nabla b\in L^2(\Di^2)$ such that
\be
\label{A-1}
\int_{\Di^2}\frac{|\nabla a|^2}{f}\ dx^2<+\infty\quad\mbox{ and }\quad\int_{\Di^2}f\,{|\nabla b|^2} dx^2<+\infty\ .
\ee
Then there exists a unique $\varphi\in W^{1,2}\cap L^\infty(\Di^2)$ satisfying
\be
\label{A-2}
\lf\{
\begin{array}{l}
\ds-\mbox{div}\lf(f\,\nabla\varphi\rg)=\nabla\, a\cdot\nabla^\perp b\ds=\p_{x_2} a\,\p_{x_1}b-\p_{x_1}a\,\p_{x_2}b\quad\mbox{ on }\Di^2\\[5mm]
\ds\varphi\ds=0\quad\mbox{ on }\p \Di^2\ .
\end{array}
\rg.
\ee
moreover the following inequality holds
\be
\label{A.3}
\|\varphi\|_{L^\infty(\Di^2)}^2+\int_{\Di^2}f\,|\nabla\varphi|^2\ dx^2\le C_0\ \int_{\Di^2}\frac{|\nabla a|^2}{f}\ dx^2\ \int_{\Di^2}f\,{|\nabla b|^2} dx^2<+\infty\ .
\ee
where $C_0>0$ is a universal constant.
\hfill $\Box$
\end{Lma}
\noindent{\bf Proof of lemma~\ref{lm-wente-weight}} First we assume $f\in C^\infty(\Di^2)$. We adopt a strategy introduced in \cite{BG} for proving Wente estimates. Let $p\in \Di^2$ and $G_p$ the Green function solution to
\be
\label{A.4}
\lf\{
\begin{array}{l}
\ds-\mbox{div}\lf(f\,\nabla G_p\rg)\ds=\, \delta_p\quad\mbox{ on }\Di^2\\[5mm]
\ds G_p\ds=0\quad\mbox{ on }\p \Di^2\ .
\end{array}
\rg.
\ee
Thanks to the maximum principle we have $G_p>0$ on $\Di^2$. For any $0\le\al<\beta$, we denote 
\be
\label{A.5}
\om_p(\al,\beta):=\lf\{x\in \Di^2\ ;\ \al\le G_p(x)\le \beta\rg\}\quad\mbox{ and }\quad \om_p(\al):=\lf\{x\in \Di^2\ ;\ \al\le G_p(x)\rg\}\ .
\ee
Since $f\in C^\infty$, thanks to classical regularity theory $G_p\in C^\infty(\Di^2\setminus \{p\})$ and using Sard's theorem
\be
\label{A.6}
\mbox{ for a. e. }\al\in {\R}_+\quad G_p^{-1}(\{\al\})\mbox{ is a regular closed curve in $\Di^2$ }.
\ee
For such an $\al>0$ we shall orient $G_p^{-1}(\{\al\})$ thanks to $$\nu:=-\frac{\nabla G_p}{|\nabla G_p|}\ .$$ Observe that we have
\be
\label{A.7}
-\int_{G_p^{-1}(\{\al\})} f\,|\nabla G_p|\ d{\mathcal H}^1=\int_{G_p^{-1}(\{\al\})} f\,\p_\nu G_p\ d{\mathcal H}^1=\int_{\om_p(\al)}\mbox{div}\lf(f\, \nabla G_p\rg)\ d{\mathcal H}^2=-1\ .
\ee
Using the coarea formula we deduce that 
\be
\label{A.8}
\mbox{for a.e. }0\le\al<\beta\quad\int_{\om_p(\al,\beta)}f\ |\nabla G_p|^2\ dx^2=\int_\al^\beta\ ds\int_{G_p^{-1}(\{s\})}  f\,|\nabla G_p|\ d{\mathcal H}^1=\beta-\al\ .
\ee
Let $\varphi$ be the $W^{1,2}\cap L^\infty$ solution to (\ref{A-2}) given in \cite{CL} (as well as in \cite{BG}). We have
\be
\label{A.9}
\varphi(p):=\int_{\Di^2}\nabla\, a\cdot\nabla^\perp b\ G_p\ dx^2
\ee
Observe that we have again thanks to the co-area formula 
\be
\label{A.10}
\begin{array}{l}
\ds\int_n^{n+1}\ ds \int_{G_p^{-1}(\{s\})} |\nabla a|\ d{\mathcal H}^1=\int_{\om_p(n,n+1)}|\nabla a|\, |\nabla G_p|\ dx^2\quad\mbox{ and }\\[5mm]
\ds\int_n^{n+1}\ ds \int_{G_p^{-1}(\{s\})} |\nabla b|\ d{\mathcal H}^1=\int_{\om_p(n,n+1)}|\nabla b|\, |\nabla G_p|\ dx^2\ .
\end{array}
\ee
Using Sard's theorem as well as the mean value formula, for any $n\in{\N}$ we can find $\al_n\in[n,n+1]$ such that $\al_n$ is a regular value for $G_p$ and
\be
\label{A.11}
\lf\{
\begin{array}{l}
\ds\int_{G_p^{-1}(\{\al_n\})} |\nabla a|\ d{\mathcal H}^1\le 2\,\int_{\om_p(n,n+1)}|\nabla a|\, |\nabla G_p|\ dx^2\\[5mm]
\ds\int_{G_p^{-1}(\{\al_n\})} |\nabla b|\ d{\mathcal H}^1\le 2\,\int_{\om_p(n,n+1)}|\nabla b|\, |\nabla G_p|\ dx^2\ .
\end{array}
\rg.
\ee
Because of the strong Hopf maximum principle we have $G_p^{-1}(\{0\})=\p \Di^2$ and we take $\al_0=0$. We write using (\ref{A.9})
\be
\label{A.12}
\begin{array}{l}
\ds\varphi(p):=\int_{\Di^2}\nabla\, a\cdot\nabla^\perp b\ G_p\ dx^2=\sum_{n\in {\N}}\int_{\om_p(\al_n,\al_{n+1})}\nabla\, a\cdot\nabla^\perp b\ G_p\ dx^2\\[5mm]
\ds\quad=\sum_{n\in {\N}}\int_{\om_p(\al_n,\al_{n+1})}\nabla\, a\cdot\nabla^\perp b\ (G_p-n)\ dx^2+\sum_{n\in {\N}}n\,\int_{\om_p(\al_n,\al_{n+1})}\nabla\, a\cdot\nabla^\perp b\ dx^2
\end{array}
\ee
Since $|G_p-n|\le 2$ on $\om_p(\al_n,\al_{n+1})$ we bound using Cauchy Schwartz
\be
\label{A.13}
\begin{array}{l}
\ds\lf|\sum_{n\in {\N}}\int_{\om_p(\al_n,\al_{n+1})}\nabla\, a\cdot\nabla^\perp b\ (G_p-n)\ dx^2\rg|\le 2\,\sum_{n\in {\N}}\int_{\om_p(\al_n,\al_{n+1})}|\nabla\, a|\,|\nabla b|\  dx^2\\[5mm]
\ds\le\, 2\, \lf[\int_{\Di^2}\frac{|\nabla a|^2}{f}\ dx^2\rg]^{1/2}\ \lf[\int_{\Di^2}f\,{|\nabla b|^2}\ dx^2\rg]^{1/2}\ .
\end{array}
\ee
Next, by integration by parts we have
\be
\label{A.14}
\begin{array}{l}
\ds\sum_{n\in {\N}}n\,\int_{\om_p(\al_n,\al_{n+1})}\nabla\, a\cdot\nabla^\perp b\ dx^2\\[5mm]
\ds=\sum_{n>0}n\,\int_{G_p^{-1}(\{\al_{n+1}\})} a\,\p_\tau b\ d\tau-\, n\,\int_{G_p^{-1}(\{\al_{n}\})} a\,\p_\tau b\ d\tau\\[5mm]
\ds=\sum_{n>1}(n-1-n)\,\int_{G_p^{-1}(\{\al_{n}\})} a\,\p_\tau b\ d\tau-\, \,\int_{G_p^{-1}(\{\al_{1}\})} a\,\p_\tau b\ d\tau\\[5mm]
\ds=-\,\sum_{n>0}\int_{G_p^{-1}(\{\al_{n}\})} a\,\p_\tau b\ d\tau\ ,
\end{array}
\ee
where $d\tau$ denotes the length element along the curve $G_p^{-1}(\{\al_{n}\})$ and $\tau$ is the positive unit tangent vector to that curve. We decompose each level set $G_p^{-1}(\{\al_{n}\})$
\[
G_p^{-1}(\{\al_{n}\})=\cup_{i\in I_n}\gamma_i^n\ ,
\]
where $\gamma_i^n$ are disjoint embedded $S^1$ in $\Di^2$ and we write
\be
\label{A.15}
\int_{G_p^{-1}(\{\al_{n}\})} a\,\p_\tau b\ d\tau=\sum_{i\in I_n} \int_{\gamma_i^n} a\,\p_\tau b\ d\tau=\sum_{i\in I_n} \int_{\gamma_i^n} (a-a_i^n)\,\p_\tau b\ d\tau\ ,
\ee
where $a_i^n$ is a value taken by $a$ on $\gamma_i^n$. Hence we have in particular
\be
\label{A.16}
\| a-a_i^n\|_{L^\infty(\gamma_i^n)}\le \int_{\gamma_i^n}|\p_\tau a|\ d\tau\ .
\ee
Hence we have using (\ref{A.8}) and (\ref{A.11})
\be
\label{A.17}
\begin{array}{l}
\ds\lf|\int_{G_p^{-1}(\{\al_{n}\})} a\,\p_\tau b\ d\tau\rg|\le \sum_{i\in I_n} \int_{\gamma_i^n}|\p_\tau a|\ d\tau\  \int_{\gamma_i^n}|\p_\tau b|\ d\tau\\[5mm]
\ds\quad\le \int_{G_p^{-1}(\{\al_{n}\})}|\nabla a|\ d{\mathcal H}^1\  \int_{G_p^{-1}(\{\al_{n}\})}|\nabla b|\ d{\mathcal H}^1\\[5mm]
\ds\quad\le 4\, \int_{\om_p(n,n+1)}|\nabla a|\, |\nabla G_p|\ dx^2\ \int_{\om_p(n,n+1)}|\nabla b|\, |\nabla G_p|\ dx^2\\[5mm]
\ds\quad\le 4 \, \int_{\om_p(n,n+1)} f\,|\nabla G_p|^2\ dx^2\ \lf[   \int_{\om_p(n,n+1)}\frac{|\nabla a|^2}{f}\, dx^2 \rg]^{\frac{1}{2}}\ \lf[   \int_{\om_p(n,n+1)}\frac{|\nabla b|^2}{f}\, dx^2 \rg]^{\frac{1}{2}}\\[5mm]
\ds\quad\le 4 \, \lf[   \int_{\om_p(n,n+1)}\frac{|\nabla a|^2}{f}\, dx^2 \rg]^{\frac{1}{2}}\ \lf[   \int_{\om_p(n,n+1)}\frac{|\nabla b|^2}{f}\, dx^2 \rg]^{\frac{1}{2}}\ .
\end{array}
\ee
Hence using Cauchy Schwartz
\be
\label{A.18}
\begin{array}{l}
\ds\lf|\sum_{n>0}\int_{G_p^{-1}(\{\al_{n}\})} a\,\p_\tau b\ d\tau\rg|\le4\,\sum_{n>0}\lf[   \int_{\om_p(n,n+1)}\frac{|\nabla a|^2}{f}\, dx^2 \rg]^{\frac{1}{2}}\ \lf[   \int_{\om_p(n,n+1)}\frac{|\nabla b|^2}{f}\, dx^2 \rg]^{\frac{1}{2}}\\[5mm]
\ds\quad\le 4\,\lf[\int_{\Di^2}\frac{|\nabla a|^2}{f}\, dx^2 \rg]^{\frac{1}{2}}\ \lf[ \int_{\Di^2}\frac{|\nabla b|^2}{f}\, dx^2 \rg]^{\frac{1}{2}}\ .
\end{array}
\ee
Combining (\ref{A.12}), (\ref{A.13}) and (\ref{A.18}) together with the fact that $f^{-1}\le f$ we finally obtain
\be
\label{A.19}
\forall p\in \Di^2\quad |\varphi(p)|\le 6\ \lf[\int_{\Di^2}\frac{|\nabla a|^2}{f}\ dx^2\rg]^{1/2}\ \lf[\int_{\Di^2}f\,{|\nabla b|^2}\ dx^2\rg]^{1/2}\ .
\ee
Multiplying the equation by $\varphi$ and integrating by parts gives, using (\ref{A.19}), 
\be
\label{A-20}
\begin{array}{l}
\ds\int_{\Di^2}f\, |\nabla\varphi|^2\ dx^2\le \|\varphi\|_{L^\infty(\Di^2)}\ \int_{\Di^2}|\nabla a|\, |\nabla b|\ dx^2\\[5mm]
\ds\quad \le 6\ \int_{\Di^2}\frac{|\nabla a|^2}{f}\ dx^2\ \int_{\Di^2}f\,{|\nabla b|^2}\ dx^2\ .
\end{array}
\ee

\medskip

Let $f\in L^\infty(\Di^2)$ and $f\ge 1$. We can extend $f$ outside $\Di^2$ by 1 and we consider $f_\ep:=\chi_\ep\star f$ where $\chi_\ep(x)=\ep^{-2}\chi(|x|/\ep)$, $\chi\in C^\infty_0(B_1)$ and $\int_{{\R}^2}\chi\ dx^2=1$. We have $a$ and $b$ in $W^{1,2}(\Di^2)$. We have clearly $f_\ep\ge 1$ and we consider $\varphi_\ep$ the solution to
\be
\label{A-21}
\lf\{
\begin{array}{l}
\ds-\mbox{div}\lf(f_\ep\,\nabla\varphi_\ep\rg)=\nabla\, a\cdot\nabla^\perp b\ds=\p_{x_2} a\,\p_{x_1}b-\p_{x_1}a\,\p_{x_2}b\quad\mbox{ on }\Di^2\\[5mm]
\ds\varphi_\ep\ds=0\quad\mbox{ on }\p \Di^2\ .
\end{array}
\rg.
\ee
Applying (\ref{A-20}) gives
\be
\label{A-22}
\begin{array}{l}
\ds\int_{\Di^2}f_\ep\, |\nabla\varphi_\ep|^2\ dx^2 \le 6\ \int_{\Di^2}\frac{|\nabla a|^2}{f_\ep}\ dx^2\ \int_{\Di^2}f_\ep\,{|\nabla b|^2}\ dx^2\ .
\end{array}
\ee
We have that $f_\ep\rightarrow f$ and $f_\ep^{-1}\rightarrow f^{-1}$ almost everywhere. Dominated convergence is then implying
\be
\label{A-23}
\lim_{\ep\rightarrow 0} \int_{\Di^2}\frac{|\nabla a|^2}{f_\ep}\ dx^2=\int_{\Di^2}\frac{|\nabla a|^2}{f}\ dx^2\quad{ and }\quad\lim_{\ep\rightarrow 0}\int_{\Di^2}f_\ep\,{|\nabla b|^2}\ dx^2=\int_{\Di^2}f\,{|\nabla b|^2}\ dx^2\ .
\ee
Hence in particular we have
\be
\label{A-24}
\begin{array}{l}
\ds\limsup_{\ep\rightarrow 0}\int_{\Di^2}f_\ep\, |\nabla\varphi_\ep|^2\ dx^2 \le 6\ \int_{\Di^2}\frac{|\nabla a|^2}{f}\ dx^2\ \int_{\Di^2}f\,{|\nabla b|^2}\ dx^2\ .
\end{array}
\ee
Modulo extraction of a subsequence we have that $\varphi_{\ep_k}\rightharpoonup \varphi$ weakly in $W^{1,2}$. Classical results on convolutions is giving that $f_\ep\rightarrow f$ strongly in $L^p(\Di^2)$ for any $p<+\infty$. Hence
\be
\label{A-25}
f_{\ep_k}\,\nabla\varphi_{\ep_k}\rightharpoonup f\,\nabla\varphi\quad\mbox{ in }{\mathcal D}'(\Di^2)\ .
\ee
We can pass to the limit in (\ref{A-21}) and we obtain
\be
\label{A-26}
\lf\{
\begin{array}{l}
\ds-\mbox{div}\lf(f\,\nabla\varphi\rg)=\nabla\, a\cdot\nabla^\perp b\ds=\p_{x_2} a\,\p_{x_1}b-\p_{x_1}a\,\p_{x_2}b\quad\mbox{ on }\Di^2\\[5mm]
\ds\varphi\ds=0\quad\mbox{ on }\p \Di^2\ .
\end{array}
\rg.
\ee
Thanks to  (\ref{A-24}) we have that $\sqrt{f_\ep} \nabla\varphi_\ep$ is uniformly bounded in $L^2$. We can then extract a subsequence from $\ep_k$ (that we keep denoting $\ep_k$) such that
\be
\label{A-27}
\sqrt{f_{\ep_k}} \nabla\varphi_{\ep_k}\ \rightharpoonup X\quad\mbox{ weakly in }{L^2}(\Di^2)\ .
\ee
Since  $f_\ep\rightarrow f$ strongly in $L^p(\Di^2)$ and since $f_\ep$ is uniformly bounded from below by 1 we write
\be
\label{A-28}
|\sqrt{f_\ep}-\sqrt{f}|\le \frac{|f_\ep-f|}{\sqrt{f_\ep}+\sqrt{f}}\le 2^{-1}\ |f_\ep-f|\ .
\ee
 We deduce that $\sqrt{f_\ep}\rightarrow \sqrt{f}$ strongly in $L^p(\Di^2)$ for any $p<+\infty$\ .
Hence
\be
\label{A-29}
\sqrt{f_{\ep_k}}\,\nabla\varphi_{\ep_k}\rightharpoonup \sqrt{f}\,\nabla\varphi\quad\mbox{ in }{\mathcal D}'(\Di^2)\ .
\ee
This implies $X=\sqrt{f}\,\nabla\varphi$ and using the lower semi-continuity of the $L^2$ norm under weak $L^2$ convergence we deduce from (\ref{A-24})
\be
\label{A-30}
\int_{\Di^2}f\, |\nabla\varphi|^2\ dx^2 \le\liminf_{k\rightarrow +\infty}\int_{\Di^2}f_{\ep_k}\, |\nabla\varphi_{\ep_k}|^2\ dx^2 \le 6\ \int_{\Di^2}\frac{|\nabla a|^2}{f}\ dx^2\ \int_{\Di^2}f\,{|\nabla b|^2}\ dx^2\ .
\ee

\medskip

Let now consider a general $f$ satisfying the hypothesis of the lemma (i.e. $f\ge 1$ and $f\in L^{p'}(\Di^2)$) and consider $\nabla a\in L^{p'}(\Di^2)$ and $\nabla b\in L^2(\Di^2)$ satisfying (\ref{A-1}).

Let $M>1$ and denote $f_M(x):=\min\{f(x),M\}$. Consider $\nu_M$ solution to
\be
\label{A-31}
\lf\{
\begin{array}{l}
\ds -\mbox{div}\lf(f_M\,\nabla\nu_M   \rg)=\mbox{div}\lf( \frac{\sqrt{f_M}}{\sqrt{f}}\,\nabla^\perp a\rg)\quad\mbox{ in }\Di^2\\[5mm]
\ds\nu_M=0\quad\mbox{ on }\p \Di^2\ .
\end{array}
\rg.
\ee
Hence, there exists $a_M$ such that
\be
\label{A-32}
 \frac{\sqrt{f_M}}{\sqrt{f}}\,\nabla^\perp a= \nabla^\perp a_M-f_M\,\nabla\nu_M\ .
\ee
Dividing by $\sqrt{f_M}$ taking the square, integrating over $\Di^2$ and using that $\nu_M=0$  on $\Di^2$ is giving
\be
\label{A-33}
\begin{array}{rl}
\ds\int_{\Di^2}\lf|\frac{\nabla a}{\sqrt{f}}\rg|^2\ dx^2&\ds=\int_{\Di^2}\lf|\frac{\nabla a_M}{\sqrt{f_M}}\rg|^2\ dx^2+\int_{\Di^2}f_M\ |\nabla \nu_M|^2\ dx^2+2\, \int_{\Di^2}\nabla^\perp a_M\cdot\nabla \nu_M dx^2\\[5mm]
 &\ds=\int_{\Di^2}\lf|\frac{\nabla a_M}{\sqrt{f_M}}\rg|^2\ dx^2+\int_{\Di^2}f_M\ |\nabla \nu_M|^2\ dx^2\ .
\end{array}
\ee
Now, taking the scalar product of (\ref{A-32}) with $\nabla\nu_M$ and integrating over $\Di^2$ is giving
\be
\label{A-34}
\begin{array}{rl}
\ds\int_{\Di^2}\frac{\sqrt{f_M}}{\sqrt{f}}\,\nabla^\perp a\cdot\nabla\nu_M\ dx^2&\ds= \int_{\Di^2}\nabla^\perp a_M\cdot\nabla\nu_M\ dx^2-\int_{\Di^2}f_M\,|\nabla\nu_M|^2\ dx^2\\[5mm]
\ds &\ds=-\int_{\Di^2}f_M\,|\nabla\nu_M|^2\ dx^2
\end{array}
\ee
From (\ref{A-33}) we deduce that $\sqrt{f_M}\,\nabla\nu_M$ is uniformly bounded in $L^2(\Di^2)$. We can extract a subsequence $M_k$ such that
\be
\label{A-35}
\sqrt{f_{M_k}}\,\nabla\nu_{M_k}\ \rightharpoonup\ Y\quad\mbox{ weakly in }L^2(\Di^2)\ .
\ee
Hence in particular
\be
\label{A-36}
\lim_{k\rightarrow +\infty}\int_{\Di^2}\frac{\sqrt{f_{M_k}}}{\sqrt{f}}\,\nabla^\perp a\cdot\nabla\nu_{M_k}\ dx^2=\int_{\Di^2} \frac{\nabla^\perp a}{\sqrt{f}}\cdot Y\ dx^2
\ee
By dominated convergence we have $f_M\rightarrow f$ strongly in $L^{p'}(\Di^2)$. Thus in particular, since $f_M\ge 1$, we have that $\sqrt{f_M}\rightarrow\sqrt{f}$ almost everywhere. By dominated convergence we deduce
\be
\label{A-37}
\sqrt{f_M}\longrightarrow \sqrt{f}\quad\mbox{ strongly in }L^{2p'}(\Di^2)\ .
\ee
Since $\nabla\nu_{M_k}$ is uniformly bounded in $L^2(\Di^2)$, there exists a subsequence (that we keep denoting $\nu_{M_k}$ and $\nu_\infty\in W^{1,2}_0(\Di^2)$ such that
\be
\label{A-38}
\nabla\nu_{M_k}\ \rightharpoonup\ \nabla\nu_\infty\quad\mbox{ weakly in }L^2(\Di^2)\ .
\ee
We deduce from the above two convergences
\be
\label{A-39}
\sqrt{f_{M_k}}\,\nabla\nu_{M_k}\ \rightharpoonup\ \sqrt{f}\,\nabla\nu_\infty\quad\mbox{ weakly in }{\mathcal D}'(\Di^2)\ .
\ee
This gives $Y=\sqrt{f}\,\nabla\nu_\infty$ and since $\nu_\infty=0$ on $\p \Di^2$ we have
\be
\label{A-40}
\int_{\Di^2} \frac{\nabla^\perp a}{\sqrt{f}}\cdot Y\ dx^2=\int_{\Di^2}\nabla^\perp a\cdot\nabla \nu_\infty\ dx^2=0\ .
\ee
Combining (\ref{A-34}), (\ref{A-36}) and (\ref{A-40}) we deduce
\be
\label{A-41}
\lim_{k\rightarrow +\infty}\int_{\Di^2}f_{M_k}\,|\nabla\nu_{M_k}|^2\ dx^2=0\ .
\ee
From (\ref{A-32}) we have
\be
\label{A-42}
\lim_{k\rightarrow +\infty}\,\int_{\Di^2}\lf| \frac{\nabla a}{\sqrt{f}} -\frac{\nabla a_{M_k}}{\sqrt{f_{M_k}}} \rg|^2\ dx^2=0\ .
\ee
Now we consider $\varphi_k$ to be the solution to
\be
\label{A-43}
\lf\{
\begin{array}{l}
\ds-\mbox{div}\lf(f_{M_k}\,\nabla\varphi_k\rg)=\nabla\, a_{M_k}\cdot\nabla^\perp b\quad\mbox{ on }\Di^2\\[5mm]
\ds\varphi_k\ds=0\quad\mbox{ on }\p \Di^2\ .
\end{array}
\rg.
\ee
Since $f_{M_k}$ is bounded in $L^\infty(\Di^2)$ for any $k$ we can apply the result for such an $f_{M_k}$ and from (\ref{A-30}) we have
\be
\label{A-44}
\int_{\Di^2}f_{M_k}\, |\nabla\varphi_k|^2\ dx^2  \le 6\ \int_{\Di^2}\frac{|\nabla a_k|^2}{f_{M_k}}\ dx^2\ \int_{\Di^2}f_{M_k}\,{|\nabla b|^2}\ dx^2\ .
\ee
By dominated convergence we have
\be
\label{A-45}
\lim_{k\rightarrow+\infty}\int_{\Di^2}f_{M_k}\,{|\nabla b|^2}\ dx^2=\int_{\Di^2}f\,{|\nabla b|^2}\ dx^2\ .
\ee
Combining (\ref{A-42}), (\ref{A-44}) and (\ref{A-45})  we obtain
\be
\label{A-46}
\limsup_{k\rightarrow +\infty}\int_{\Di^2}f_{M_k}\, |\nabla\varphi_k|^2\ dx^2  \le 6\ \int_{\Di^2}\frac{|\nabla a|^2}{f}\ dx^2\ \int_{\Di^2}f\,{|\nabla b|^2}\ dx^2\ .
\ee
 We deduce that $\sqrt{f_{M_k}}\,\nabla\varphi_k$ is uniformly bounded in $L^2(\Di^2)$. We can extract a subsequence still denoted $M_k$ such that
\be
\label{A-47}
\sqrt{f_{M_k}}\,\nabla\varphi_k\ \rightharpoonup\ Z\quad\mbox{ weakly in }L^2(\Di^2)\ .
\ee
Since $\nabla\varphi_{k}$ is uniformly bounded in $L^2(\Di^2)$, there exists a subsequence (that we keep denoting $\varphi_{k}$ and $\varphi\in W^{1,2}_0(\Di^2)$ such that
\be
\label{A-48}
\nabla\varphi_{k}\ \rightharpoonup\ \nabla\varphi\quad\mbox{ weakly in }L^2(\Di^2)\ .
\ee
Using the convergence (\ref{A-37}) we deduce
\be
\label{A-49}
\sqrt{f_{M_k}}\,\nabla\varphi_k\ \rightharpoonup\ \sqrt{ f}\, \nabla \varphi\quad\mbox{ weakly in }{\mathcal D}'(\Di^2)\ .
\ee
Hence we deduce
\be
\label{A-50}
Z=\sqrt{ f}\, \nabla \varphi\ .
\ee
Hence from (\ref{A-47}) we have
\be
\label{A-51}
\sqrt{f_{M_k}}\,\nabla\varphi_k\ \rightharpoonup\ \sqrt{ f}\, \nabla \varphi\ \quad\mbox{ weakly in }L^2(\Di^2)\ .
\ee
The lower semicontinuity of the $L^2$ norm under weak convergence is implying
\be
\label{A-52}
\int_{\Di^2}f\, |\nabla\varphi|^2\ dx^2 \limsup_{k\rightarrow +\infty}\int_{\Di^2}f_{M_k}\, |\nabla\varphi_k|^2\ dx^2  \le 6\ \int_{\Di^2}\frac{|\nabla a|^2}{f}\ dx^2\ \int_{\Di^2}f\,{|\nabla b|^2}\ dx^2\ .
\ee
Combining now the weak $L^2$ convergence of $\sqrt{f_{M_k}}\,\nabla\varphi_k$ towards $\sqrt{ f}\, \nabla \varphi$ with the strong convergence of $\sqrt{f_{M_k}}$ towards $\sqrt{f}$ in $L^2$ we deduce
\be
\label{A-53}
f_{M_k}\,\nabla\varphi_k=\sqrt{f_{M_k}}\,\sqrt{f_{M_k}}\,\nabla\varphi_k\ \rightharpoonup\ \sqrt{f}\,\sqrt{f}\,\nabla\varphi ={ f}\, \nabla \varphi\quad\mbox{ weakly in }{\mathcal D}'(\Di^2)\ .
\ee
This implies that $\varphi$ is a solution of 
\be
\label{A-54}
\lf\{
\begin{array}{l}
\ds-\mbox{div}\lf(f\,\nabla\varphi\rg)=\nabla\, a\cdot\nabla^\perp b\quad\mbox{ on }\Di^2\\[5mm]
\ds\varphi\ds=0\quad\mbox{ on }\p \Di^2\ .
\end{array}
\rg.
\ee
and $\varphi$ satisfies the inequality (\ref{A-46}). This concludes the proof of lemma~\ref{lm-wente-weight}.\hfill $\Box$
\section{Some estimates for $A $-Wente type systems}
\reset
\begin{Lma}
\label{lm-wente-type-estimates} There exists $\sigma>0$ depending only on $m\in {\N}$ such that, for any  $A\in L^{\infty}\cap W^{1,2}(\Di^2,M_m({\R}))$ satisfying
\be
\label{B-1}
\lf\|\mbox{dist}(A,SO(m))\rg\|^2_{L^\infty(\Di^2)}+\int_{\Di^2}|\nabla A|^2\ dx^2\le \sigma
\ee
then, for any $p\in[2,3]$ and any $D\in W^{1,p'}(\Di^2,M_m({\R})$ and any $v\in W^{1,p}(\Di^2,{\R}^m)$ there exists a unique $\varphi\in W^{1,2}(\Di^2,{\R}^m)$ such that 
\be
\label{B-2}
\lf\{
\begin{array}{l}
\ds-\mbox{div}(A\,\nabla\varphi)=\nabla^\perp D\cdot\nabla v\quad\mbox{ in }{\mathcal D}'(\Di^2)\\[5mm]
\ds \varphi=0\quad\mbox{on }\p \Di^2\ .
\end{array}
\rg.
\ee
moreover there exists a constant $C>0$ independent of $p\in [2,3]$ such that
\be
\label{B-2-a}
\|\nabla\varphi\|_{L^{2,1}(\Di^2)}\le C\, \|\nabla D\|_{L^{p'}(\Di^2)}\ \|\nabla v\|_{L^{p}(\Di^2)}\  .
\ee
\hfill $\Box$
\end{Lma}
\noindent{\bf Proof of lemma~\ref{lm-wente-type-estimates}.} First we establish the existence and uniqueness of $\varphi$ in $W^{1,q}_0(\Di^2,{\R}^m)$ for any $q<2$. The system
(\ref{B-2}) is equivalent to
\be
\label{B-3}
\lf\{
\begin{array}{l}
\ds-\Delta(A\varphi)=-\mbox{div}(\nabla A\, A^{-1}\, A\varphi)+\nabla^\perp D\cdot\nabla v\quad\mbox{ in }{\mathcal D}'(\Di^2)\\[5mm]
\ds A\varphi=0\quad\mbox{on }\p \Di^2\ .
\end{array}
\rg.
\ee
We fix $q<2$ and we construct inductively $w_k$ such that
\be
\label{B-4}
\lf\{
\begin{array}{l}
\ds-\Delta w_0=\nabla^\perp D\cdot\nabla v\quad\mbox{ in }{\mathcal D}'(\Di^2)\\[5mm]
\ds w_0=0\quad\mbox{on }\p \Di^2\ .
\end{array}
\rg.
\ee
and for $k>0$
\be
\label{B-5}
\lf\{
\begin{array}{l}
\ds-\Delta w_k=-\mbox{div}(\nabla A\, A^{-1}\, w_{k-1})\quad\mbox{ in }{\mathcal D}'(\Di^2)\\[5mm]
\ds w_k=0\quad\mbox{on }\p \Di^2\ .
\end{array}
\rg.
\ee
We have, since $w_0=0$ on $\p \Di^2$
\be
\label{B-6}
\|\nabla w_0\|_{L^q(\Di^2)}\le C\,\|\nabla w_0\|_{L^{2,\infty}(\Di^2)}\le C\, \|\Delta w_0\|_{L^1(\Di^2)}\le C\, \|\nabla D\|_{L^{p'}(\Di^2)}\ \|\nabla v\|_{L^p(\Di^2)}
\ee
And
\be
\label{B-7}
\|\nabla w_k\|_{L^q(\Di^2)}\le C_q\, \|\nabla A\|_{L^2(\Di^2)}\, \|w_{k-1}\|_{L^{q^\ast}(\Di^2)}\ ,
\ee
where $(q^\ast)^{-1}+2^{-1}=q^{-1}$. Sobolev embedding is giving
\be
\label{B-8}
\|w_{k-1}\|_{L^{q^\ast}(\Di^2)}\le C^\ast_q\ \|\nabla w_{k-1}\|_{L^q(\Di^2)}\ .
\ee
Hence we obtain
\be
\label{B-9}
\|\nabla w_k\|_{L^q(\Di^2)}\le C'_q\, \|\nabla A\|_{L^2(\Di^2)}\,\ \|\nabla w_{k-1}\|_{L^q(\Di^2)}\ .
\ee
For $C'_q\, \|\nabla A\|_{L^2(\Di^2)}<1$ the sum $w:=\sum w_k$ is converging in $W^{1,q}_0(\Di^2)$ and $\varphi= A^{-1} w$ is a solution to the system.
We proceed to the following Hodge decomposition in $L^q(\Di^2)$
\be
\label{B-10}
A\,\nabla \varphi=\nabla w_0+\nabla^\perp\zeta\ .
\ee
Using the above estimate we prove that $\varphi$ solution to (\ref{B-2}) is unique in $W^{1,q}_0(\Di^2)$.

\medskip

 Recall that from Coifman Lions Meyer and Semmes \cite{CLMS} (see also \cite{Hel} section 3.3) for any choice of $\nabla a\in L^{3/2}(\Di^2)$ and $\nabla b\in L^3(\Di^2)$, the solution $\phi$ to
\be
\label{B-11}
\lf\{
\begin{array}{rl}
\ds -\Delta\phi&\ds=\nabla^\perp a\cdot\nabla b\quad\mbox{ in }\Di^2\\[5mm]
\ds\phi&\ds=0\quad\mbox{ on }\p \Di^2\ .
\end{array}
\rg.
\ee
Statisfies
\be
\label{B-12}
\|\nabla\phi\|_{L^{2,1}(\Di^2)}\le \|\nabla^2\phi\|_{L^1(\Di^2)}\le C\ \|\nabla a\|_{L^{3/2}(\Di^2)}\  \|\nabla b\|_{L^{3}(\Di^2)}
\ee
Similarly, exchanging the role of $\nabla a$ and $\nabla b$ gives
\be
\label{B-13}
\|\nabla\phi\|_{L^{2,1}(\Di^2)}\le \|\nabla^2\phi\|_{L^1(\Di^2)}\le C\ \|\nabla a\|_{L^{3}(\Di^2)}\  \|\nabla b\|_{L^{3/2}(\Di^2)}
\ee
The bilinear interpolation theorem (section 4.4  together with theorem 5.3.1 of \cite{BL}  ) gives that there exists a constant $C>0$ independent of $p\in[3/2,3]$ such that
\be
\label{B-14}
\|\nabla\phi\|_{L^{2,1}(\Di^2)}\le C\ \|\nabla a\|_{L^{p}(\Di^2)}\  \|\nabla b\|_{L^{p'}(\Di^2)}
\ee
Applying this inequality to $w_0$ is giving
\be
\label{B-15}
\|\nabla w_0\|_{L^{2,1}(\Di^2)}\le C\ \|\nabla D\|_{L^{p'}(\Di^2)}\ \|\nabla v\|_{L^{p}(\Di^2)}\  .
\ee
The map $\zeta$ satisfies
\be
\label{B-16}
\lf\{
\begin{array}{l}
\ds-\Delta \zeta=\nabla^\perp A\cdot\nabla \varphi\quad\mbox{ in }{\mathcal D}'(\Di^2)\\[5mm]
\ds \frac{\p\zeta}{\p\nu}=0\quad\mbox{on }\p \Di^2\ .
\end{array}
\rg.
\ee
A-priori, there is no Wente type estimate (i.e $L^\infty\cap W^{1,2}$) for the Neumann boundary conditions (see \cite{DP}) but we are in the very special case where
$\varphi=0$ on $\p \Di^2$ and then Wente estimates hold true (see theorem 1.3 of \cite{DP}). Hence we have the a-priori estimate
\be
\label{B-17}
\begin{array}{l}
\ds\|\nabla\zeta\|_{L^{2,1}(\Di^2)}\le C\, \|\nabla A\|_{L^2(\Di^2)}\ \|\nabla\varphi\|_{L^2(\Di^2)}\\[5mm]
\ds\quad\le C\, \|\nabla A\|_{L^2(\Di^2)}\ \|\nabla\zeta\|_{L^2(\Di^2)}+C\, \|\nabla A\|_{L^2(\Di^2)}\ \|\nabla w_0\|_{L^2(\Di^2)}
\end{array}
\ee
Hence for $C\, \|\nabla A\|_{L^2(\Di^2)}<1$ we obtain the a-priori estimate
\be
\label{B-18}
\|\nabla\zeta\|_{L^{2,1}(\Di^2)}\le C\, \|\nabla A\|_{L^2(\Di^2)}\ \|\nabla w_0\|_{L^2(\Di^2)}\ .
\ee
Combined with (\ref{B-15}) this gives the a-priori estimate
\be
\label{B-19}
\|\nabla\varphi\|_{L^{2,1}(\Di^2)}\le C\, \|\nabla B\|_{L^{p'}(\Di^2)}\ \|\nabla v\|_{L^{p}(\Di^2)}\  .
\ee
The result follows by density.\hfill $\Box$
\begin{Lma}
\label{lm-wente-type-estimates-bis} There exists $\sigma>0$ depending only on $m\in {\N}$ such that, for any  $A\in L^{\infty}\cap W^{1,2}(\Di^2,M_m({\R}))$ satisfying
\be
\label{B-20}
\lf\|\mbox{dist}(A,SO(m))\rg\|^2_{L^\infty(\Di^2)}+\int_{\Di^2}|\nabla A|^2\ dx^2\le \sigma
\ee
then, for any $p\in[2,3]$ and any $D\in W^{1,2}(\Di^2,M_m({\R}))$ and any $v\in W^{1,{(2,\infty)}}(\Di^2,{\R}^m)$ there exists a unique $\varphi\in W^{1,2}(\Di^2,{\R}^m)$ such that 
\be
\label{B-21}
\lf\{
\begin{array}{l}
\ds-\mbox{div}(A\,\nabla\varphi)=\nabla^\perp D\cdot\nabla v\quad\mbox{ in }{\mathcal D}'(\Di^2)\\[5mm]
\ds \varphi=0\quad\mbox{on }\p \Di^2\ .
\end{array}
\rg.
\ee
moreover there exists a constant $C>0$ independent of $p\in [2,3]$ such that
\be
\label{B-22}
\|\nabla\varphi\|_{L^{2}(\Di^2)}\le C\, \|\nabla D\|_{L^{2}(\Di^2)}\ \|\nabla v\|_{L^{2,\infty}(\Di^2)}\  .
\ee
\hfill $\Box$
\end{Lma}
\noindent{\bf Proof of Lemma~\ref{lm-wente-type-estimates-bis}.} First we establish the existence and uniqueness of $\varphi$ in $W^{1,q}_0(\Di^2,{\R}^m)$ for any $q<2$. The system
(\ref{B-21}) is equivalent to
\be
\label{B-23}
\lf\{
\begin{array}{l}
\ds-\Delta(A\varphi)=-\mbox{div}(\nabla A\, A^{-1}\, A\varphi)+\nabla^\perp D\cdot\nabla v\quad\mbox{ in }{\mathcal D}'(\Di^2)\\[5mm]
\ds A\varphi=0\quad\mbox{on }\p \Di^2\ .
\end{array}
\rg.
\ee
We fix $q<2$ and we construct inductively $w_k$ such that
\be
\label{B-24}
\lf\{
\begin{array}{l}
\ds-\Delta w_0=\nabla^\perp D\cdot\nabla v\quad\mbox{ in }{\mathcal D}'(\Di^2)\\[5mm]
\ds w_0=0\quad\mbox{on }\p \Di^2\ ,
\end{array}
\rg.
\ee
and for $k>0$
\be
\label{B-25}
\lf\{
\begin{array}{l}
\ds-\Delta w_k=-\mbox{div}(\nabla A\, A^{-1}\, w_{k-1})\quad\mbox{ in }{\mathcal D}'(\Di^2)\\[5mm]
\ds w_k=0\quad\mbox{on }\p \Di^2\ .
\end{array}
\rg.
\ee
We have, since $w_0=0$ on $\p \Di^2$
\be
\label{B-26}
\|\nabla w_0\|_{L^q(\Di^2)}\le C\,\|\nabla w_0\|_{L^{2}(\Di^2)}\le  C\, \|\nabla D\|_{L^{2}(\Di^2)}\ \|\nabla v\|_{L^{(2,\infty}(\Di^2)}\ ,
\ee
and
\be
\label{B-27}
\|\nabla w_k\|_{L^q(\Di^2)}\le C_q\, \|\nabla A\|_{L^2(\Di^2)}\, \|w_{k-1}\|_{L^{q^\ast}(\Di^2)}\ ,
\ee
where $(q^\ast)^{-1}+2^{-1}=q^{-1}$. Sobolev embedding is giving
\be
\label{B-28}
\|w_{k-1}\|_{L^{q^\ast}(\Di^2)}\le C^\ast_q\ \|\nabla w_{k-1}\|_{L^q(\Di^2)}\ .
\ee
Hence we obtain
\be
\label{B-29}
\|\nabla w_k\|_{L^q(\Di^2)}\le C'_q\, \|\nabla A\|_{L^2(\Di^2)}\,\ \|\nabla w_{k-1}\|_{L^q(\Di^2)}\ .
\ee
For $C'_q\, \|\nabla A\|_{L^2(\Di^2)}<1$ the sum $w:=\sum w_k$ is converging in $W^{1,q}_0(\Di^2)$ and $\varphi= A^{-1} w$ is a solution to the system.
We proceed to the following Hodge decomposition in $L^q(\Di^2)$
\be
\label{B-30}
A\,\nabla \varphi=\nabla w_0+\nabla^\perp\zeta\ .
\ee
Using the above estimate we prove that $\varphi$ solution to (\ref{B-2}) is unique in $W^{1,q}_0(\Di^2)$. The map $\zeta$ satisfies
\be
\label{B-31}
\lf\{
\begin{array}{l}
\ds-\Delta \zeta=\nabla^\perp A\cdot\nabla \varphi\quad\mbox{ in }{\mathcal D}'(\Di^2)\\[5mm]
\ds \frac{\p\zeta}{\p\nu}=0\quad\mbox{on }\p \Di^2\ .
\end{array}
\rg.
\ee
and using  (\ref{B-18}), for $C\, \|\nabla A\|_{L^2(\Di^2)}<1$ we obtain the a-priori estimate
\be
\label{B-32}
\|\nabla\zeta\|_{L^{2,1}(\Di^2)}\le C\, \|\nabla A\|_{L^2(\Di^2)}\ \|\nabla w_0\|_{L^2(\Di^2)}\ .
\ee
Combining (\ref{B-26}) and (\ref{B-32}) we get the result. \hfill $\Box$

\section{Explicit Constant in H\"older-Lorentz inequality }
\reset
We recall a result below which is well known to the experts in it's qualitative formulation but for which we need to have the explicit constant.
\begin{Lma}
\label{lm-holder-lorentz}
Let $p>2$. There exists a constant $C$ universal such that, for any $g\in L^p(\Di^2)$, there holds
\be
\label{A-54-a}
\|g\|_{L^{2,1}(\Di^2)}\le C \ \pi^{\frac{p-2}{2p}}\ \lf( \frac{2\,(p-1)}{p-2}\rg)^{\frac{p-1}{p}}\ \|g\|_{L^p(\Di^2)}\ .
\ee
\hfill $\Box$
\end{Lma}
\noindent{\bf Proof of Lemma~\ref{lm-holder-lorentz}.} Let $p>2$ and $g\in L^p(\Di^2)$. We have
\be
\label{A-55}
\|g\|_{L^p(\Di^2)}^p=p\, \int_0^{+\infty}t^{p-1}\, \lf|\lf\{x\in \Di^2\ ;\ |g(x)|>t\rg\}\rg|\ dt\ .   
\ee
We denote by $g^\ast(t)$\footnote{Given $f\colon\Di^2\to \R$ measurable,  we define its decreasing rearrangement
$f^{\ast} \colon [0,+\infty]\to  [0,+\infty]$ as
$$f^\ast(t):=\inf\{0\le\lambda\le +\infty:~~ \lf|\lf\{x\in \Di^2\ ;\ |f(x)|>\lambda\rg\}\rg|\le t\}$$
with the convention that $0\cdot\infty=\infty\cdot 0=0$ (as it is customary in measure theory). }
the decreasing rearrangement of $|g|$ such that forall $t\in {\R}_+$
\be
\label{A-56}
\lf|\lf\{x\in \Di^2\ ;\ |g(x)|>t\rg\}\rg|=\lf|\lf\{s\in [0,\pi]\ ;\ g^\ast(s)>t\rg\}\rg|\ .
\ee
We recall
\be
\label{A-57}
\|g\|_{L^{2,1}(\Di^2)}\simeq \int_0^\pi g^\ast(t)\ \frac{dt}{\sqrt{t}}\ ,
\ee
and from (\ref{A-55}) and (\ref{A-56}) we get
\be
\label{A-58}
\|g\|_{L^p(\Di^2)}^p= p\,\int_0^{+\infty}t^{p-1}\, \lf|\lf\{s\in [0,\pi]\ ;\ |g^\ast(s)|>t\rg\}\rg|\ dt=\int_0^\pi\, ( g^\ast)^p(t)\ dt\ .
\ee
This gives in particular
\be
\label{A-59}
\begin{array}{l}
\ds\|g\|_{L^{2,1}(\Di^2)}\le C\, \int_0^\pi g^\ast(t)\ \frac{dt}{\sqrt{t}} \\[5mm]
\ds\quad\le C\, \lf( \int_0^\pi \frac{dt}{t^{\frac{p}{2(p-1)}}}  \rg)^{\frac{p-1}{p}}\ \lf(\int_0^\pi\, ( g^\ast)^p(t)\ dt\rg)^{\frac{1}{p}} \\[5mm]
\ds\quad\le C\, \lf( \frac{2\,(p-1)}{p-2}\ \pi^{\frac{p-2}{2(p-1)}}  \rg)^{\frac{p-1}{p}}\ \|g\|_{L^p(\Di^2)}\,.
\end{array}
\ee
This proves the lemma.\hfill $\Box$
 
\end{document}